%% file: jpurcell-cusps.tex
\documentclass[11pt]{amsart}

\usepackage{amsfonts, amstext, amsmath, amsthm, amscd, amssymb}
\usepackage{graphicx, color, epsfig, psfrag}

\newcommand{\tild}[1]{{\widetilde{#1}}}

\newcommand{\HH}{{\mathbb{H}}}
\newcommand{\RR}{{\mathbb{R}}}

\newcommand{\CC}{{\mathbb{C}}}

\def\co{\colon\thinspace}

\theoremstyle{plain}
\newtheorem{theorem}{Theorem}[section]
\newtheorem{corollary}[theorem]{Corollary}
\newtheorem{lemma}[theorem]{Lemma}
\newtheorem{prop}[theorem]{Proposition}

\newtheorem*{no-num-theorem}{Theorem}

\theoremstyle{definition}
\newtheorem{define}[theorem]{Definition}

\newtheorem*{remark}{Remark}
\newtheorem*{notation}{Notation}

\begin{document}

\title{Cusp Shapes Under Cone Deformation}
\author{Jessica S. Purcell}
\thanks{This research was supported in part by NSF grant DMS-0704359.}

\address{Department of Mathematics, Brigham Young
	University, Provo, UT 84602, U.S.A}
\address{Mathematical Institute, University of Oxford, 24-29 St
	Giles', Oxford, OX1 3LB, U.K.}
\email{jpurcell@math.byu.edu}

\begin{abstract}
  \input{abstract}

\end{abstract}

\maketitle


\section{Introduction \label{sec:intro}}
\input{introduction}

\section{Background on Cone Deformation \label{sec:deform-cone}}
\input{deform-cone}

\section{Change in Cusp Shape \label{sec:cusp-est}}
\input{cusp-est}

\section{Bounding the boundary terms \label{sec:bdrybound}}
\input{bdrybound.tex}

\section{Bounding the meridian \label{sec:merbound}}
\input{meridian.tex}

\section{Cusp shapes of hyperbolic knots \label{sec:knots}}
\input{hyp-knot}

\providecommand{\bysame}{\leavevmode\hbox to3em{\hrulefill}\thinspace}
\providecommand{\href}[2]{#2}

\end{document}

%% file: abstract.tex
A horospherical torus about a cusp of a hyperbolic manifold inherits a
Euclidean similarity structure, called a cusp shape.  We bound the
change in cusp shape when the hyperbolic structure of the manifold is
deformed via cone deformation preserving the cusp.  The bounds are in
terms of the change in structure in a neighborhood of the singular
locus alone.

We then apply this result to provide information on the cusp shape of
many hyperbolic knots, given only a diagram of the knot.
Specifically, we show there is a universal constant C such that if a
knot admits a prime, twist reduced diagram with at least C crossings
per twist region, then the length of the second shortest curve on the
cusp torus is bounded.  The bound is linear in the number
of twist regions of the diagram.

%% file: introduction.tex

\subsection{Motivation}


The interplay between geometry and topology is an important theme in
the study of 3--manifolds.  Since Thurston's work in the late 1970's,
it has been known that any Haken 3--manifold has a geometric
decomposition \cite{thurston}.  Recently Perelman has announced a
proof that any closed, orientable 3--manifold has such a decomposition
(the geometrization conjecture).  However, it is important not only to
know the existence of a geometric structure on a 3--manifold, but also
to understand the relation between that structure and the topological
description of the manifold.

Often a 3--manifold is described topologically, for example by a
combinatorial diagram indicating a knot complement in $S^3$, or a Dehn
filling description (see Definition \ref{def:dehn} below).  However,
relating this topological description to the geometry of the manifold
seems to be very difficult, and remains an important problem in the
area.  There are few tools available to give geometric data from
topological information.  In this paper, we present a new such tool.

We are particularly interested in topological descriptions of
3--manifolds with torus boundary.  These include all knot and link
complements in $S^3$, as well as knot and link complements in other
3--manifolds.  It is a classic result due to Wallace and Lickorish
that every closed, orientable 3--manifold is given by a Dehn filling
on a link complement in $S^3$ (\cite{wallace}, \cite{lickorish}).
Thus these manifolds are of fundamental importance in the subject.

We will also focus only on those manifolds which are hyperbolic.  The
link complements of Wallace and Lickorish above can be taken to be
hyperbolic.  Thurston showed that on any such cusped hyperbolic
manifold ``most'' Dehn fillings are hyperbolic.  That is, only a
finite number of slopes per link component need be excluded, then any
Dehn filling on remaining slopes yields a hyperbolic manifold
\cite{thurston}.  Additionally, for knots in $S^3$, any knot that is
not a torus or satellite knot is known to be hyperbolic
\cite{thurston:bulletin}.  Therefore in some sense, most 3--manifolds
are hyperbolic.  Mostow--Prasad rigidity implies that the hyperbolic
structure on these finite volume manifolds is unique \cite{mostow},
\cite{prasad}.  However, this unique structure is little understood.

Given a hyperbolic manifold with torus boundary components, i.e. given
a \emph{cusped} hyperbolic manifold, the geometric structure on the
manifold allows us to associate a Euclidean similarity structure,
called a cusp shape, to each boundary component.  The cusp shapes of a
manifold are geometric quantities which we will relate to a manifold's
topological description.

There are several reasons one might be interested in knowing the cusp
shape of a manifold.  Cusp shape is an interesting piece of geometric
information by itself.  It also can be used to give a rough lower
bound on the volume of the manifold, by calculating cusp volumes and
using packing arguments such as those of B{\"o}r{\"o}czky
\cite{boroczky}.  (See however \cite{eudave-munoz-luecke}: this lower
bound may not be very sharp.)  More importantly, knowing the cusp
shape, we also obtain information on the geometry of manifolds
obtained by Dehn filling.  For example, we can find all filling slopes
which might not yield a hyperbolic manifold \cite{agol:bounds},
\cite{lackenby:surg}, \cite{hk:univ}.

Our main results give explicit bounds on the change in cusp shape of a
3--manifold with multiple cusps under all but finitely many Dehn
fillings, where one cusp is left unfilled.  This is essentially the
content of Theorems \ref{thm:main-cusp}, \ref{thm:main-norm-length},
and \ref{thm:main-norm-height} below.  The first, Theorem
\ref{thm:main-cusp}, is the most general, bounding cusp information.
In Theorem \ref{thm:main-norm-length}, we give a differential
inequality bounding the normalized lengths of curves on cusps.  In
Theorem \ref{thm:main-norm-height}, we bound the change in height of
the cusp shape.

We also apply these results to the case of knots in $S^3$.  Given only
a diagram of a large class of knots, we are able to determine bounds
on the height of cusps.  This is the content of Theorem
\ref{thm:main-knot-cusp}.  Note that similar results, relating the
geometry of a knot to its diagram, are rare.  In fact, we know only of
results relating volumes to a diagram.  Lackenby found lower bounds on
volumes of alternating knots based on a condition of a diagram, and
upper bounds for general knots \cite{lackenby:alt-volume}.  The lower
bounds were improved by Agol, Storm and Thurston \cite{ast:volumes},
and the upper by Agol and D. Thurston in an appendix to Lackenby's
paper \cite{lackenby:alt-volume}.  Using cone deformations, we found
similar bounds on volumes for another class of knots
\cite{purcell:volume}.  Using negatively curved metrics, we recently
improved these results with Futer and Kalfagianni \cite{fkp:volumes}.

However, as far as we are aware, the results of this paper are the
first of their kind giving actual geometric bounds on cusp shape based
only on a diagram.

\subsection{Cone deformations}

Our primary technique is cone deformation, which we will now describe.
Using properties of cone deformations, we obtain geometric information
on the topological process of Dehn filling.

\begin{define}
  Given a manifold $M$ with torus boundary $\partial M$ and specified
  slope $s$ on $\partial M$, we obtain a new manifold by gluing a
  solid torus to $\partial M$ such that $s$ bounds a disk in the solid
  torus.  This is called the \emph{Dehn filling of $M$ along $s$}, or
  sometimes \emph{Dehn surgery}, and the resulting manifold is denoted
  $M(s)$.  More generally, if $M$ has $n$ torus boundary components
  with slopes $s_1, \dots, s_n$, then we obtain $M(s_1, \dots, s_n)$
  in the same manner.
\label{def:dehn}
\end{define}

A geometric description of Dehn filling was investigated by Thurston
\cite{thurston}.  Given a complete hyperbolic structure on a manifold
$M$ with torus boundary, he showed the structure could be deformed
through incomplete hyperbolic structures with singularities at the
core of the added solid torus to structures whose completion is again
non-singular.  These structures are exactly hyperbolic structures on
Dehn filled manifolds $M(s)$.  The space of deformations is called
hyperbolic Dehn surgery space.

We restrict our attention to particular paths through hyperbolic Dehn
surgery space in which the incomplete hyperbolic manifolds are all
cone manifolds.  That is, the completion gives a manifold with cone
singularities along the core of the added solid torus, with cone angle
$\alpha$.  (A more complete definition of a cone deformation is given
in Definition \ref{def:cone-deform} below.)

We have some control over cone deformations.  Local rigidity results
due to Hodgson and Kerckhoff \cite{hk:cone-rigid}, extended by
Bromberg in the case of geometrically finite manifolds
\cite{bromberg:cone-rigid}, give control on the change in geometry of
the entire manifold, knowing only information on the change in a
neighborhood of the singular locus.

This geometric control has been used recently to prove several
interesting results in 3--manifold theory.  Hodgson and Kerckhoff used
it to give universal bounds on the number of exceptional Dehn surgery
slopes for finite volume manifolds \cite{hk:univ}.  In the infinite
volume case, cone deformations lie at the heart of drilling and
grafting theorems, described by Bromberg in \cite{bromberg:cone}.
These have recently been used to resolve problems in the field of
Kleinian groups, such as the Bers density conjecture
\cite{brock-bromberg:density}, and existence of type preserving
sequences of Kleinian groups \cite{brock-bromberg-evans-souto:tame}.
This paper gives another application of cone deformations, to finite
volume manifolds with cusps.

\subsection{Cusp shapes}

Let $M$ be a manifold admitting a hyperbolic structure, with a torus
boundary component $T$. Under the hyperbolic structure on $M$, $T$
will be a cusp.  We take a horoball neighborhood of this cusp.  Its
boundary is a horospherical torus, which inherits a Euclidean
structure from the hyperbolic structure on $M$.  The structure is
independent of choice of horoball neighborhood, up to scaling of the
Euclidean metric.  Thus we obtain a class of Euclidean similarity
structures on $T$, which we call the \emph{cusp shape} of the torus
$T$.

Now, suppose the hyperbolic structure of $M$ is changing under a cone
deformation with singularities along a link $\Sigma$ in $M$ disjoint
from $T$.  $\Sigma$ is called the singular locus of the deformation.
$T$ will remain a cusp throughout the deformation, and its cusp shape
will be changing with the change of geometry of $M$.  We bound this
change in terms of the change in geometry in a neighborhood of
$\Sigma$.

\begin{theorem}  Suppose we have the following:
  \begin{itemize}
  \item $M$ a complete finite volume hyperbolic 3--manifold with a
    cusp $T$.
  \item $X=X_0$ a complete finite volume hyperbolic 3--manifold which
    can be joined to $M=X_\tau$ by a smooth one--parameter family of
    hyperbolic cone manifolds $X_t$, such that cone angles go from $0$
    at time $t=0$ to $2\pi$ at time $\tau$.
  \item A singular locus, say with $n$ components, such that the tube
    radius about each component is larger than $R_1\geq 0.56$ for all
    $t$.
  \item A horoball neighborhood $U_t$ of the cusp of $X_t$ which
    deforms to $T$ at time $\tau$.  Let $\gamma(t)$ denote the path
    $\partial U_t$ travels through the Teichm{\"u}ller space of the
    torus, endowed with Teichm{\"u}ller metric $\| \cdot \|$.
  \end{itemize}
	Then there exists a parameterization of the deformation such that
  we have the following bound for $\partial U_t$:
  $$2\,\|\gamma'(t)\|^2 \mbox{Area}(\partial U_t) \leq n\,C(R_1)$$
  Here $C(R_1)$ is a strictly decreasing function of $R_1$ which
  approaches $0$ as $R_1$ approaches infinity.
  Moreover, we can assume $\tau \leq (2\pi)^2$.
\label{thm:main-cusp}
\end{theorem}

Also interesting is the application of Theorem \ref{thm:main-cusp}
to lengths of curves on the cusp in question.  Let $\beta$ be a slope
on the cusp torus, and again let $\gamma(t)$ be the path that the cusp
shape $\partial U_t$ traces through the {T}eichm{\"u}ller space of the
torus.  Let $L_t$ be the normalized length of $\beta$ under the
deformation, with derivative $\dot{L} = \frac{d}{dt}L_t$ (see
section \ref{sec:normlength} for more information).  Then $L_t$ will
satisfy:
\begin{equation}
  -\|\gamma'(t)\|\, L_t \leq \dot{L} \leq
   \|\gamma'(t)\|\, L_t.
  \label{eq:norm-length-change}
\end{equation}

Combining equation (\ref{eq:norm-length-change}) with the result of
Theorem \ref{thm:main-cusp}, we may bound the change in normalized
length of a curve under a cone deformation.

\begin{theorem}
Let the set up be as in Theorem \ref{thm:main-cusp} above.  Let $L_t$
be the normalized length of a slope on a cusp torus under the cone
deformation, with derivative $\dot{L} = \frac{d}{dt}L_t$.  Then
there exists a parameterization of the cone deformation such that
the change in $L_t$ is bounded by the following inequality:
$$ -\left(\sqrt{\frac{n\,C(R_1)}{2\mbox{Area}(\partial U_t)}}\right)
\, L_t \leq \dot{L} \leq
\left(\sqrt{\frac{n\,C(R_1)}{2\mbox{Area}(\partial U_t)}}\right) \,
L_t
$$
\label{thm:main-norm-length}
\end{theorem}

Finally, we use Theorem \ref{thm:main-norm-length} to bound the
normalized length of the arc perpendicular to the shortest curve on
the cusp torus, and then bound the length of this shortest curve
below.  We then obtain an estimate that is independent of time.

\begin{theorem}
  Let the notation be as in Theorem \ref{thm:main-cusp} above.  Denote
  by $h(M)$ the normalized \emph{height} of the cusp torus $\partial
  U_\tau$.  That is, $h(M)$ is the normalized length of the arc
  perpendicular to the shortest curve on $\partial U_\tau$ that runs
  from the shortest curve back to the shortest curve (this arc might
  not be closed).  Similarly, let $h(X)$ denote the normalized height
  of $\partial U_0$.  Then
	there exists a parameterization of the cone deformation such that
	the change in normalized height is bounded in terms of $R_1$ alone:
$$-(2\pi)^2\frac{\sqrt{n\,C(R_1)}}{(1-e^{-2R_1})\sqrt{2}} \leq
 h(M) - h(X) \leq
(2\pi)^2\frac{\sqrt{n\,C(R_1)}}{(1-e^{-2R_1})\sqrt{2}}.$$
\label{thm:main-norm-height}
\end{theorem}




\subsection{Applications to knot theory}

In the last part of this paper we give specific applications of
Theorem \ref{thm:main-norm-height}.  We determine bounds on the cusp
shape of a knot based solely on a diagram of that knot.

Before we state the results, we review some definitions.

\begin{define}
  Let $K \subset S^3$ be a knot with diagram $D(K)$.  Define a
  \emph{twist region} to be a region of $D(K)$ where two strands of
  the diagram twist around each other maximally.  Precisely, if $D(K)$
  is viewed as a 4--valent planar graph with over--under crossing
  information at each vertex, then a twist region is a maximal string
  of bigons in the complement of the graph arranged end to end, or a
  single vertex of the graph adjacent to no bigons.  
\end{define}

We also require the diagram to be reduced, in the sense of the
following two definitions.  These definitions are illustrated
elsewhere (see e.g. \cite{lackenby:surg}, \cite{futer-purcell}).

\begin{define}
  A diagram $D(K)$ is \emph{prime} if when any simple closed curve
  $\gamma$ intersects $D(K)$ transversely in two points in the
  interiors of edges, then $\gamma$ bounds a subdiagram containing no
  crossings of the original diagram.
	\label{def:prime}
\end{define}

\begin{define}
  The diagram $D(K)$ is \emph{twist reduced} if when any simple closed
  curve $\gamma$ intersects $D(K)$ transversely in four points in the
  interiors of edges, with two of these points adjacent to one
  crossing and two others adjacent to another crossing, then $\gamma$
  bounds a subdiagram consisting of a (possibly empty) collection of
  bigons arranged end to end between these two crossings.
	\label{def:twred}
\end{define}

Any diagram of a prime knot or link can be simplified into a prime,
twist reduced diagram.  If the diagram is not prime, then crossings on
one side of a simple closed curve are extraneous and can be removed.
If the diagram is not twist reduced, then performing flypes will
amalgamate two twist regions adjacent to a simple closed curve into
one, reducing the number of twist regions.

For a hyperbolic knot complement, we can choose a horoball
neighborhood about the cusp such that a meridian on the boundary of
that horoball has length $1$ (see e.g. \cite{adams:waist}).  Consider
the geodesic arc orthogonal to the meridian which runs from meridian
to meridian on this horospherical torus.  Define the length of this
arc to be the \emph{height} of the cusp.  Note that when the meridian
is the shortest curve on the cusp (as will be the case in the knots of
this paper), this definition of height agrees with that in Theorem
\ref{thm:main-norm-height}.

Let $K$ be a knot in $S^3$ which admits a prime, twist reduced diagram
with at least $2$ twist regions, with each twist region containing at
least $6$ crossings.  Then it was shown in \cite{futer-purcell} that
$S^3-K$ is hyperbolic.  The following theorem gives bounds on the
cusp shape.

\begin{theorem}
  Let $K$ be a knot in $S^3$ which admits a prime, twist reduced
  diagram with a total of $n\geq 2$ twist regions, with each twist
  region containing at least $c\geq 116$ crossings.  In a hyperbolic
  structure on $S^3-K$, take a horoball neighborhood $U$ about $K$.
  Normalize so that the meridian on $\partial U$ has length $1$.  Then
  the height $H$ of the cusp of $S^3-K$ satisfies:
  $$H\geq n\,(1-f(c))^2.$$
  Here $f(c)$ is a positive function of $c$ which approaches $0$ as
  $c$ increases to infinity.

  Additionally, $H\leq n(\sqrt{n-1} + f(c))^2$.
\label{thm:main-knot-cusp}
\end{theorem}

The function $f$ is given explicitly in Section \ref{sec:knots}.  We
can plug in values of $c$ to obtain more specific estimates.  For
example, in Corollary \ref{cor:0.8} we let $c=145$.  Then $1-f(c)$ is
approximately $0.8154$.

Note that the shortest non-meridional slope on the cusp of $S^3-K$
will have length at least that of the height.  So Theorem
\ref{thm:main-knot-cusp} implies that in this case, the shortest
non-meridional slope has length at least $n\,(1-f(c))^2$.

By results of Hodgson and Kerckhoff \cite{hk:univ}, Dehn filling along
a slope with normalized length at least $7.515$ results in a hyperbolic
manifold.  Corollary \ref{cor:0.8} implies that for a complicated
knot, the shortest non-meridional slope will have normalized length at
least $\sqrt{n}(0.8154)$.  Thus the normalized length of a slope will
be larger than $7.515$ whenever $n\geq 85$.  Thus we have the
following corollary to Theorem \ref{thm:main-knot-cusp}:

\begin{corollary}
  Let $K$ be a knot in $S^3$ which admits a prime, twist reduced
  diagram with at least $85$ twist regions and at least $145$
  crossings per twist region.  Then any non-trivial {D}ehn filling of
  $S^3-K$ is hyperbolic.
\label{cor:knot-cusp}
\end{corollary}

Corollary \ref{cor:knot-cusp} is proven without any reference to the
geometrization conjecture.  If we assume that conjecture, then we
recently proved (with Futer \cite{futer-purcell}) that the numbers
$85$ and $145$ could be reduced to $4$ and $6$ respectively.  However,
note the proof of that result gives no geometric information on cusp
geometry.

The numbers of crossings, 116 and 145, can be reduced significantly
for certain classes of knots.  We illustrate this in section
\ref{subsec:2bridge} for 2--bridge knots.

\subsection{Organization of this paper}
\label{sec:organize}

In Section \ref{sec:deform-cone}, we review background information on
cone deformations.

In Section \ref{sec:cusp-est}, we begin the proofs of Theorems
\ref{thm:main-cusp}, \ref{thm:main-norm-length}, and
\ref{thm:main-norm-height} using information set up in Section
\ref{sec:deform-cone}.  This leads us to results which contain certain
boundary terms.

In Section \ref{sec:bdrybound}, we give bounds on these boundary
terms, which allow us to conclude the proof of Theorems
\ref{thm:main-cusp} and \ref{thm:main-norm-length}.

We complete the proof of \ref{thm:main-norm-height} in Section
\ref{sec:merbound}.

Finally, in Section \ref{sec:knots} we apply the results of the
previous sections to knots, and present the proof of Theorem
\ref{thm:main-knot-cusp}.

\subsection{Acknowledgments}

We would like to thank Steve Kerckhoff for many helpful conversations.
We lean heavily on his results with Craig Hodgson.  His willingness to
clarify, explain, and correct has made this paper possible.  We would
also like to thank the referee for very carefully reading this paper,
and for many helpful suggestions, corrections, and comments.


%% file: deform-cone.tex
In this section, we will recall definitions associated with cone
deformations and review relevant results from work of Hodgson and
Kerckhoff (\cite{hk:cone-rigid}, \cite{hk:univ}).

\begin{define}
Let $M$ be a 3--manifold containing a closed 1--manifold $\Sigma$.  A
\emph{hyperbolic cone structure} on $M$ is an incomplete hyperbolic
structure on $X=M-\Sigma$, whose completion is singular along
$\Sigma$.  A cross section of a component of $\Sigma$ is a hyperbolic
cone with angle $\alpha$, where $\alpha$ is constant along that
component.  (For a more complete description of the metric on $M$, see
\cite{hk:cone-rigid}.)  The set $\Sigma$ is called the \emph{singular
locus} of $M$.  $M$ with this structure is called a \emph{hyperbolic
cone manifold}.  A \emph{hyperbolic cone deformation} on $M$ is a
smooth one--parameter family of hyperbolic cone structures $X_t$.
\label{def:cone-deform}
\end{define}

\subsection{Infinitesimal deformations}
Now, we may study deformations by considering possible ``derivatives
at time $t$'' of $X_t$, or infinitesimal deformations of hyperbolic
structure,
for any time $t$.
Each of these infinitesimal deformations corresponds to an element in
a certain cohomology group $H^1(X; E)$, where recall $X=M-\Sigma$.  We
will be manipulating elements of $H^1(X;E)$, so we recall some
information on their structure.


\subsubsection{Cohomology}

Given $\omega$ in $H^1(X;E)$, $\omega$ is a one--form on $X$ with
values in the vector bundle $E$ of infinitesimal isometries of $X$.
The bundle $E$ has fiber $sl(2,\CC)$, which is the {L}ie algebra of
infinitesimal isometries of $\HH^3$.  Recall that this {L}ie algebra
has a complex structure.  Geometrically, if $s$ represents an
infinitesimal translation in the direction of $s$, then $i s$
represents an infinitesimal rotation with axis in the direction of
$s$.  Thus on $X$ we can identify $E$ with the complexified tangent
bundle $TX \otimes \CC$, and write any $\omega \in H^1(X;E)$ in
complex form: $\omega = v + iw$.

By the {H}odge theorem proved in \cite{hk:cone-rigid}, in each
cohomology class of $H^1(X;E)$ there is a \emph{harmonic}
representative of the form
$$\omega = \eta + i*D\eta,$$
where $\eta$ is a unique $TX$-valued
1--form on $X = M-\Sigma$ that satisfies:
$$D^*\eta=0,$$
\begin{equation}
D^*D\eta + \eta =0.
\label{eqn:d*deta}
\end{equation}
Here $D$ is the exterior covariant derivative on such forms and $D^*$
is its adjoint.

Thus we may assume that our cone deformation is chosen such that at
each time $t$, the infinitesimal deformation of hyperbolic structure
of $X_t$ corresponds to a \emph{harmonic} element $\omega$ in
$H^1(X;E)$.

\subsubsection{Tubular and Horoball neighborhoods}

For purposes of this paper, we are particularly interested in the
effect of the deformation in a neighborhood of the singular locus or
in a horoball neighborhood of a cusp.  So in this subsection, we will
review a decomposition of $\omega$ (and $\eta$) into a nice form in
these types of neighborhoods.

Select a component of the singular locus $\Sigma$ of $M$.  Let $V$ be
a tubular neighborhood of that component, with tube radius $R$.  $V$
is a solid torus with a singular core curve.  During the cone
deformation, the geometric structure on $V$ is changing in meridional,
longitudinal, and radial directions.  We can write $\omega$ in $V$ to
reflect this:
\begin{equation}
\omega = \omega_0 + \omega_c.
\label{eqn:omega-decomp}
\end{equation}
Here only $\omega_0 = \eta_0 + i*D\eta_0$
changes the holonomy of the meridian and longitude on the torus
$\partial V$, and $\omega_c$ is a correction term.

We can further decompose $\omega_0$ as follows:
\begin{equation}
\omega_0 = u\, \omega_m + (x+iy)\,\omega_\ell.
\label{eqn:omega_0-decomp}
\end{equation}
Here, $u$, $x$, and $y$ are real numbers.  The 1--forms $\omega_m$ and
$\omega_\ell$ are standard forms, calculated in \cite{hk:cone-rigid},
which depend only on the tube radius $R$ of $V$.  The form $\omega_m =
\eta_m + i*D\eta_m$ gives the change in the meridional direction.  In
particular, it represents the infinitesimal deformation which
decreases the cone angle but doesn't change the real part of the
complex length of the meridian.  The form $\omega_\ell = \eta_\ell
+i*D\eta_\ell$ stretches the singular locus, but leaves the holonomy
of the meridian unchanged.

Since we will use the explicit formula for the standard form
$\omega_\ell$ in a calculation in a later section, we will restate
that form here.  Let $e_1$, $e_2$, and $e_3$ be an orthonormal frame
for $V$ in cylindrical coordinates, with $e_1$ pointing in the radial
direction and $e_2$ tangent to the meridian.  Let $R$ be the radius of
the tube $V$.  Recall $\omega_\ell$ is a ($TX\otimes\CC$)--valued
1--form, which can be viewed as an element of $Hom(TX, TX\otimes\CC)$.
Thus it can be described as a matrix in the cylindrical coordinates
$e_1$, $e_2$, $e_3$.
\begin{equation}
  \omega_\ell(R) = \left[
    \begin{array}{ccccc}
      (\cosh R)^{-2} & \hspace{.1in} & 0 &\hspace{.1in}& 0 \\
      0 && -1 && -i\tanh R \\
      0 && -i\tanh R && 1+ (\cosh R)^{-2}
    \end{array}
    \right]
  \label{eqn:omega_l}
\end{equation}

Now, so far, our decomposition of $\omega$ has been in a tubular
neighborhood $V$ of a component of the singular locus.  We can do a
similar decomposition of $\omega$ in a horoball neighborhood $U$ of a
cusp.  In fact, we may consider $U$ as a tubular neighborhood of a
1--dimensional submanifold of $M$ on which the cone angle remains $0$
throughout the deformation, and the tube radius is infinite.  Again we
may decompose $\omega$ as $\omega_0 + \omega_c$, where only $\omega_c$
affects the holonomy of the torus $\partial U$.  But now, since the
cone angle at the core of $U$ is not changing, the meridional part of
$\omega_0$ is identically $0$.  Thus in this case,
$$\omega_0 = (x + iy)\, \omega_\ell,$$
where again $x$ and $y$ are real numbers, and $\omega_\ell$ is a
standard form.  Since the tube radius of $U$ is infinite, to write
down the formula for $\omega_\ell$ explicitly we let $R$ go to
infinity in equation (\ref{eqn:omega_l}).
\begin{equation}
  \omega_\ell(\infty) = \left[
    \begin{array}{ccc}
      0 & 0 & 0\\
      0 & -1 & -i\\
      0 & -i & 1
    \end{array}\right]
  \label{eqn:omega_l-cusp}
\end{equation}


\subsection{Boundary terms}
\label{subsec:boundary-terms}

Our proof of Theorem \ref{thm:main-cusp} will involve certain boundary
terms.  In this section, we define these terms.  We explain how they
arise and recall certain properties.  They will be used in the
proof of Theorem \ref{thm:main-cusp}.


Using results of the last section, we may always assume that we have a
harmonic representative $\omega = \eta + i*D\eta$.  Thus $\eta$
satisfies equation (\ref{eqn:d*deta}).

Let $N$ be any submanifold of $X = M-\Sigma$ with boundary $\partial
N$ oriented by the outward normal.  When we take the $L^2$ inner
product of equation (\ref{eqn:d*deta}) with $\eta$ and integrate by
parts over $N$, we obtain the boundary formula (\cite[Lemma
2.3]{hk:univ}):
\begin{equation}
  B(N)= \int_{\partial N} \eta \wedge *D\eta = \| \eta \|^2_N + \|
  *D\eta \|^2_N 
  \label{eqn:BN}
\end{equation}
Here $\| \cdot \|$ denotes the $L^2$--norm on $N$.  We will apply
equation (\ref{eqn:BN}) to submanifolds of $X=M-\Sigma$ involving
tubular or horoball neighborhoods.

First, we will introduce some notation.  Let $V_j$ be a tubular
neighborhood of a component (the $j$-th component) of the singular
locus, or a horoball neighborhood of a cusp.  Using similar notation
to that of \cite{hk:univ}, define
$$b_{V_j}(\zeta,\xi)=\int_{\partial V_j}*D\zeta\wedge \xi,$$
where
$\partial V_j$ is oriented such that the tangent vector orthogonal to
$\partial V_j$ is an outward normal when viewed from $V_j$.  We will
be considering $b_{V_j}(\eta, \eta)$ in this paper, and we review
results concerning this term.

The decomposition (\ref{eqn:omega-decomp}) of $\omega$ in $V_j$ into
$\omega_0 + \omega_c$, and $\eta$ into $\eta_0 + \eta_c$ decomposes
the term $b_{V_j}(\eta, \eta)$:
\begin{equation}
b_{V_j}(\eta, \eta) = b_{V_j}(\eta_0, \eta_0) +
b_{V_j}(\eta_c, \eta_c).
\label{eqn:bV-decomp}
\end{equation}
That is, cross terms vanish (\cite[Lemma 2.5]{hk:univ}).
Additionally, the term $b_{V_j}(\eta_c, \eta_c)$ is actually
non-positive (\cite[Lemma 2.6]{hk:univ}), so we find
\begin{equation}
  b_{V_j}(\eta, \eta) \leq b_{V_j}(\eta_0, \eta_0)
\label{eqn:eta<eta_0}
\end{equation}

Now, let $V$ be a tubular neighborhood of the entire singular locus
$\Sigma$, as well as horoball neighborhoods of cusps.  Thus $V$ is a
union of tubular neighborhoods $V_j$ of components $\Sigma_j$ of the
singular locus and horoball neighborhoods $U_j$.  Then letting $N=X-V$
in equation (\ref{eqn:BN}),
\begin{eqnarray}
B(X-V)  & = &\| \eta \|^2_{X-V} + \|*D\eta\|^2_{X-V} =
 \int_{\partial V} *D\eta \wedge \eta \\
 & = &
 \sum_j b_{V_j} (\eta, \eta) + \sum_k b_{U_j}(\eta,\eta).
\label{eqn:sumbj=eta}
\end{eqnarray}
The first sum is over all components of the singular locus, the second
over all horoball neighborhoods of cusps.  This equation, along with
equation (\ref{eqn:eta<eta_0}), implies
\begin{equation}
  0 \leq  \sum_j b_{V_j} (\eta, \eta)+ \sum_k
  b_{U_k}(\eta,\eta)  \leq \sum_j b_{V_j} (\eta_0,
  \eta_0) + \sum_k b_{U_k} (\eta_0, \eta_0).
\label{eqn:sum>0}
\end{equation}

Notice $\eta_0 = (\eta_0)_j$ depends on the neighborhood $V_j$ or
$U_j$.  However, using the notation $b_{V_j}(\eta_0, \eta_0)$, it
should be clear that we are referring to the decomposition of $\eta$
particular to $V_j$ in this context, and similarly for $U_j$.  We will
further simplify notation, using the following definition.

\begin{define}
  Let $\Sigma_j$ be a component of the singular locus $\Sigma$.  Let
  $V_j$ be a tubular neighborhood of $\Sigma_j$ of radius $R$.  We let
  $b_j$ be the boundary term:
  $$b_j = b_{V_j}(\eta_0, \eta_0).$$
  \label{def:bj}
\end{define}

\begin{remark}
Note that for our applications of the results of this section, the
one--form $\omega$ corresponds to an infinitesimal deformation of
hyperbolic structure at time $t$.  Thus $\omega$, $\eta$, $\eta_0$,
$b_j$, etc. will all depend on time $t$ in our applications in future
sections.

To avoid a notational nightmare, in the sequel we will assume the time
$t$ has been fixed, and continue writing these terms without reference
to $t$, except occasionally where we feel recalling the dependency on
$t$ will help avoid confusion.
\end{remark}

%% file: cusp-est.tex
In Section \ref{sec:deform-cone}, we set up notation and reviewed
known results.  Given this information, we are ready to begin the
proof of Theorem \ref{thm:main-cusp}.

\subsection{Boundary relations}

Restrict to the case when we have a single cusp.

We will let $U$ be a horoball neighborhood of the one cusp, and $V_1$,
\dots, $V_n$ be tubular neighborhoods of the components of the
singular locus $\Sigma$, each with radius $R$.  Here, we are letting
$n$ be the total number of components of $\Sigma$.  We will assume
throughout that these neighborhoods are chosen such that all
intersections of distinct neighborhoods are trivial.

Let $\omega \in H^1(X;E)$ correspond to a harmonic infinitesimal
deformation of hyperbolic structure at time $t$, and decompose
$\omega$ in tubular neighborhoods and horoball neighborhoods as in the
previous section, and consider boundary terms.

By equation (\ref{eqn:sum>0}), we see immediately that
\begin{equation}
-b_U(\eta_0, \eta_0) \leq \sum_{j=1}^n b_j.
\label{eqn:eta-bound}
\end{equation}
Again, note that $\eta_0$ in the term $b_U(\eta_0,\eta_0)$ refers to
the decomposition of $\eta$ valid only in the horoball neighborhood
$U$.
Also, recall that the terms of equation (\ref{eqn:eta-bound}) all
depend on time $t$, but we are suppressing $t$ for notational
purposes.

Now, we will analyze the left hand side of (\ref{eqn:eta-bound}).
Recall that in a horoball neighborhood of a cusp, $\omega_0$ can be
written as a complex multiple of the form $\omega_\ell =
\omega_\ell(\infty)$ in equation (\ref{eqn:omega_l-cusp}).  That is,
there are real numbers $a$ and $b$
(for fixed $t$)
such that
\begin{equation}
\omega_0 = (a+ib)\left[
  \begin{array}{ccc}
    0&0&0\\
    0&-1&-i\\
    0&-i&1
  \end{array}\right]
\label{eqn:eta_0-matrix}
\end{equation}

\begin{lemma}
  Using the notation above,
  $$-b_U(\eta_0, \eta_0) = 2(a^2+b^2)\mbox{Area}(\partial U).$$
  \label{lemma:b_U=area}
\end{lemma}

\begin{proof}
  By equation (\ref{eqn:eta_0-matrix}), we may write $\eta_0 =
  \mbox{Re}(\omega_0)$ and $*D\eta_0 = \mbox{Im}(\omega_0)$ as:
\begin{equation}
  \eta_0 = \left[\begin{array}{ccc}
      0&0&0\\
      0&-a&b\\
      0&b&a\end{array}\right] =
  \left[\begin{array}{c}0\\-a\\b\end{array}\right] \omega_2 +
  \left[\begin{array}{c}0\\b\\a\end{array}\right] \omega_3
\label{eqn:eta0mat}
\end{equation}
and:
$$
  *D\eta_0 = \left[\begin{array}{ccc}
      0&0&0\\
      0&-b&-a\\
      0&-a&b\end{array}\right] =
  \left[\begin{array}{c}0\\-b\\a\end{array}\right] \omega_2 +
  \left[\begin{array}{c}0\\-a\\b\end{array}\right] \omega_3.
$$
Here $\omega_1$, $\omega_2$, and $\omega_3$ are forms dual to the
vectors $e_1$, $e_2$, and $e_3$ giving cylindrical coordinates on
$U$.  Recall in particular that $e_1$ is radial, and $e_2$ is
tangent to a meridian.

We then may compute $-b_U(\eta_0, \eta_0)$ explicitly in terms of the
constants $a$ and $b$.
$$-b_U(\eta_0, \eta_0) = \int_{\partial U} \eta_0\wedge *D\eta_0 =
\int_{\partial U} 2(a^2+b^2)\omega_2\wedge\omega_3.$$

Since $\omega_2\wedge\omega_3$ is the area form for the torus
$\partial U$, we have:
$$-b_U(\eta_0, \eta_0) = 2 (a^2+b^2)\mbox{Area}(\partial U).$$
\end{proof}

Notice that equation
(\ref{eqn:eta-bound}) and Lemma \ref{lemma:b_U=area} imply:
\begin{equation}
0 \leq \sum_j b_j,
\label{eqn:b_j>0}
\end{equation}
and that we have equality by equation (\ref{eqn:sumbj=eta}) if and
only if $\eta$ is trivial on $X-V$, which happens only when the
deformation is trivial.

\subsection{Teichm{\"u}ller Space}
\input{teichmuller}

%% file: teichmuller.tex
We now turn our attention to the term $(a^2+b^2)$ in Lemma
\ref{lemma:b_U=area}.  This term is closely related to the change
in metric of the torus $\partial U$.

Under the hyperbolic metric on $X=M-\Sigma$, the boundary $\partial U$
inherits a Euclidean structure.  This gives a point in the
Teichm\"uller space of the torus, $\mathcal{T}(T^2)$.  

\begin{remark}
We view $\mathcal{T}(T^2)$ as the space of flat structures on the
torus.  To each flat structure $\zeta$ is associated a unit area
Euclidean metric $g_{\zeta}$, and vice versa.  However, in the
following discussion we will refer to both the flat structure and its
associated metric, and so we will continue to distinguish between the
two.  We use Greek letters $\gamma$, $\zeta$ to denote points in
$\mathcal{T}(T^2)$, and $g$ with appropriate subscripts to denote the
induced metric on the torus.
\end{remark}

As the metric on $X$ changes under the cone deformation, the point in
$\mathcal{T}(T^2)$ will also change, tracing out a smooth path
$\gamma(t)$, with tangent vector $\gamma'(t)$.  The tangent vector
$\gamma'(t)$ is an infinitesimal change of the flat structure on the
torus.  Our 1--form $\omega$ encodes the infinitesimal deformation of
the hyperbolic cone structure $X_t$ on $X=M-\Sigma$.  So we may relate
$\gamma'(t)$ to $\omega$, or more particularly, to $a^2+b^2$.  We will
do so using the Teichm\"uller metric on $\mathcal{T}(T^2)$.

Definitions of the Teichm\"uller distance vary by constant factors in
the literature.  We will use the convention that the Teichm\"uller
distance between two points, $\sigma$ and $\tau$ in
$\mathcal{T}(T^2)$, is defined to be
$$d(\sigma, \tau) = \frac{1}{2}\inf_{f}\;\log K_f,$$ where $f\co
\sigma \to \tau$ is a $K_f$--quasiconformal map, with $K_f$ the
smallest such constant.  (See for example \cite{lehto}.)


\begin{lemma}
  Let $\gamma(t)$ denote the path of $\partial U$ through
  $\mathcal{T}(T^2)$.  Let $\gamma'(t)$ denote its tangent vector.
  Then
  $$\sqrt{a^2+b^2} = \left\| \gamma'(t) \right\|,$$
  where the metric is the {T}eichm{\"u}ller metric.
\label{lemma:a+b-gamma}
\end{lemma}

\begin{proof}
Fix $t=t_0$, and consider the flat structure on the torus given by
$\partial U_{t_0} = \gamma(t_0)$.  We may assume it has unit area.
Recall that any Teichm\"uller geodesic through $\gamma(t_0)$ is given
as follows.  Select two orthogonal geodesic foliations $F_1$ and $F_2$
on the torus.  Select $\lambda>0$.  For any $u \in \RR$, obtain a new
metric by multiplying vectors along $F_1$ by $\exp(-\lambda\,u)$ and
multiplying vectors along $F_2$ by $\exp(\lambda\,u)$.  This is a
``stretch--squeeze'' map, which preserves the area of the torus.  It
gives a one--parameter family of flat structures on the torus, which
we denote by $\zeta(u)$.

Teichm\"uller's theorem says that in the case of the torus, this
stretch--squeeze map has the minimal quasiconformal distortion, so is
a Teichm\"uller geodesic.  For each $u$, the derivative of the map
$\zeta(u)$ takes an infinitesimal circle to an infinitesimal ellipse,
the ratio of whose axes is $\exp(2u\lambda)$.  Thus for any $u$,
the distance between the structure $\zeta(0) = \gamma(t_0)$ at time
$0$ and $\zeta(u)$ at time $u$ is
$$d(\zeta(0), \zeta(u)) = \frac{1}{2}\inf_f\log K_f =
\frac{1}{2}\log(e^{2u\lambda}) = \lambda\,u.$$ 
Therefore, $\lambda = \| \zeta'(u) \|$, where the norm is given with
respect to the Teichm\"uller metric.  In particular, $\lambda =
\|\zeta'(0)\|$.

On the other hand, for any point on the torus, let $u_1$ and $u_2$ be
orthonormal vectors in the directions of $F_1$ and $F_2$ respectively.
Then when we write the infinitesimal change of metric induced by
$\zeta'(0)$ in coordinates given by $u_1$ and $u_2$, we obtain at
every point a diagonal matrix $\Lambda$ with $-\lambda$, $\lambda$ on
the diagonal, since at each time $u$ those vectors are just multiplied
by $e^{-\lambda\,u}$, $e^{\lambda\,u}$, respectively.  Thus if we let
$g_u$ be the metric on the torus induced by $\zeta(u)$, we have
$$\left.\frac{d}{du}\right|_{u=0} g_u(x,y) = 2 g_0(\Lambda\,x, y).$$

Now we relate this discussion to the term $a^2+b^2$.

Recall the infinitesimal deformation of the hyperbolic cone structure
of $X$ is encoded by the one--form $\omega$.  In particular, the
infinitesimal change in metric is given by the real part $\eta$.  That
is, if we let $g_t$ be the metric on $X$ at time $t$, then 
\begin{equation}
\left.\frac{d}{dt}\right|_{t_0}\, g_t(x,y) \;=\; 2 g_{t_0}(\eta(t_0) x, y).
\label{eqn:dotg}
\end{equation}
(See the displayed equation on page 374 of \cite{hk:univ}.)

We are interested in the infinitesimal change in the metric on
$\partial U_{t_0}$.  Since $\eta$ decomposes into $\eta_0 + \eta_c$,
and only $\eta_0$ changes the holonomy of $\partial U$, we know the
infinitesimal change in metric on $\partial U_{t_0}$ is given by
$\eta_0$.

We have an explicit formula for $\eta_0$, written in orthonormal
coordinates with respect to the metric $g_{t_0}$, given by $e_1$,
$e_2$, $e_3$, where $e_1$ is radial and $e_2$ points in the direction
of the meridian.  Consider again equation (\ref{eqn:eta0mat}).  Since
the first row and column for the matrix of $\eta_0$ are zero, we see
that $\eta_0$ actually has no effect on or contribution to the change
of metric in the radial direction of $U$.  Hence $\eta_0$ itself gives
the infinitesimal change of metric on the torus $\partial U_{t_0}$, or
$\gamma'(t_0) = \eta_0(t_0)$.  

When we write $\eta_0$ in orthonormal coordinates $e_2$, $e_3$ on the
torus, we discard the first row and column of the matrix of equation
(\ref{eqn:eta0mat}), and obtain the $2\times 2$ matrix
\begin{equation}
A = \left[\begin{array}{cc}-a&b\\b&a\end{array}\right].
\label{eqn:define-A}
\end{equation}
Here the
first column of $A$ gives the infinitesimal change in the meridional
direction, the second the change in the direction orthogonal to the
meridian.  That is, if we restrict the metric $g_t$ to
$\partial U_t$, we may rewrite equation (\ref{eqn:dotg}) as
$$\left.\frac{d}{dt}\right|_{t=t_0}\, g_t(x,y) =
2g_{t_0}(\eta_0(t_0)\,x, y) = 2\left<Ax, y\right>.$$

Hence $A$ is a matrix representation of $\gamma'(t_0)$, the
infinitesimal change of metric of the torus $\partial U_{t_0}$.

We can diagonalize $A$.  For any point on the torus, there exist
orthonormal vectors $u_1$ and $u_2$ such that when we put $A$ into
coordinates given by $u_1$ and $u_2$, $A$ is diagonal with
$-\sqrt{a^2+b^2}$ and $\sqrt{a^2+b^2}$ on the diagonal.  Moreover,
since $A$ does not depend on the point on the torus (that is, for any
point on the torus we obtain the same matrix $A$), $u_1$ and $u_2$
determine orthogonal geodesic foliations $F_1$ and $F_2$.  These, in
turn, determine a Teichm\"uller geodesic $\zeta(u)$ with $\zeta(0) =
\gamma(t_0)$, whose initial tangent vector $\zeta'(0)$ agrees with $A$
at every point, and has Teichm\"uller norm $\|\zeta'(0\| =
\sqrt{a^2+b^2}$.  Thus infinitesimally, $\gamma(t)$ and $\zeta(u)$
agree near $t=t_0$ and $u=0$, respectively.  Hence $\|\gamma'(t_0)\| =
\|\zeta'(0)\| = \sqrt{a^2+b^2}$.
\end{proof}

\begin{remark}
In the proof of Lemma \ref{lemma:a+b-gamma}, we showed that for any
time $t$, $\gamma(t)$ agrees infinitesimally with an explicit
Teichm\"uller geodesic which is completely determined by the matrix
$A$.  We will use this again below.
\end{remark}

Finally, putting Lemmas \ref{lemma:b_U=area}, and
\ref{lemma:a+b-gamma} together with equation (\ref{eqn:eta-bound}), we
have completed the proof of the following theorem.

\begin{theorem}
	Let $X_t$ be a hyperbolic cone deformation with a cusp.  That is,
  for each time $t$, the hyperbolic cone manifold $X_t$ has a cusp
  which remains a cusp throughout the deformation.  Let $U_t$ be a
  horoball neighborhood of the cusp.  Let $\gamma(t)$ denote the path
  $\partial U_t$ travels through the Teichm{\"u}ller space of the
  torus.  Then the change in $\partial U_t$ is bounded by the following
  inequality.
$$2\,\|\gamma'(t)\|^2 \mbox{Area}(\partial U_t) \leq \sum b_j(t)$$
The sum on the right hand side is over all components of the singular
locus.  
\label{thm:cusp-1}
\end{theorem}

Note Theorem \ref{thm:cusp-1} gives an inequality identical to that of
Theorem \ref{thm:main-cusp}, except for the boundary terms $b_j$ on
the right hand side.  We will bound these, and complete the proof of
Theorem \ref{thm:main-cusp}, in Section \ref{sec:bdrybound}.  First,
we set up a similar theorem to Theorem \ref{thm:main-norm-length}.

\subsection{Normalized lengths on the torus}
	\label{sec:normlength}

We are interested in determining how \emph{normalized lengths} of
curves on the torus $\partial U$ change over the deformation. 

\begin{define}
  Let $\beta$ be a slope on the torus $T$, that is, an isotopy class
  of simple closed curves.  The normalized length of $\beta$ is
  defined to be:
  $$\mbox{Normalized length}(\beta) =
  \frac{\mbox{Length}(\beta)}{\sqrt{\mbox{Area}(T)}},$$
	where $\mbox{Length}(\beta)$ is defined to be the length of a
  geodesic representative of $\beta$ on $T$.
\end{define}

\begin{lemma}
	Let $\beta$ be a slope on $\partial U_t$.  Let $L_t(\beta)$ denote
	its normalized length, with derivative $\dot{L}(\beta)$. Then
	$L_t(\beta)$ satisfies:
\begin{equation}
  -\|\gamma'(t)\| \, L_t(\beta) \;\leq\; \dot{L}(\beta) \;\leq\;
   \|\gamma'(t)\| \,L_t(\beta).
  \label{eqn:norm-length-L}
\end{equation}
\label{lemma:ldot-bound}
\end{lemma}

\begin{proof}
Fix a time $t=t_0$.  We will show the lemma for $t_0$.

We saw in the proof of Lemma \ref{lemma:a+b-gamma} that $\gamma(t)$
agrees infinitesimally with an explicit Teichm\"uller geodesic
$\zeta(u)$ determined by the matrix $A$ of equation
(\ref{eqn:define-A}) when $t=t_0$ and $u=0$.  More specifically,
$\zeta(u)$ multiplies vectors along $F_1$ by
$\exp(-u\sqrt{a^2+b^2})$, and multiplies vectors along $F_2$ by
$\exp(u\sqrt{a^2+b^2})$, where $F_1$ and $F_2$ are orthogonal
geodesic foliations (determined by the matrix $A$) on the torus
$\partial U_{t_0}$.  Thus we will prove the lemma by showing it is
true when the change in metric on $\partial U_{t_0}$ is given by
$\zeta(u)$.

Now, $L_t(\beta)$ is defined to be the length of a geodesic
representative of $\beta$ divided by $\sqrt{\rm{Area}(\partial U_t)}$.
Without loss of generality, we may assume that $\rm{Area}(\partial
U_t) = 1$, since $\gamma'(t)$ is given by a trace free (hence area
preserving) matrix $A$.  Thus $L_t(\beta)$ is just the length of a
geodesic representative of $\beta$ at time $t_0$.

Let $\hat{\beta}_{t_0}\co I \to \partial U_{t_0}$ be a geodesic
representative of $\beta$ at time $t_0$.  Then $\hat{\beta}_{t_0}$
makes some angle $\theta$ with the foliation $F_1$.  For any $u$,
$\zeta(u)$ takes $\hat{\beta}_{t_0}$ to a new geodesic, still with
slope $\beta$, and length $L_u(\beta) = \alpha_u(\theta) \,
L_{t_0}(\beta)$.  Here
$$\alpha_u(\theta)
= \sqrt{\cos^2\theta \,\exp(-2u\sqrt{a^2+b^2}) +\sin^2\theta \exp(2u\sqrt{a^2+b^2})}$$
is maximized when $\theta=\pi/2$ and minimized when $\theta=0$.  Thus
we have
$$\exp(-\sqrt{a^2+b^2}\,u)\,L_{t_0}(\beta) \leq L_u(\beta) \leq
\exp(\sqrt{a^2+b^2}\,u) \,L_{t_0}(\beta).$$

Since $\left.\frac{d}{dt}\right|_{t=t_0} L_t(\beta) =
\left.\frac{d}{du}\right|_{u=0} L_u(\beta),$
$$-\sqrt{a^2+b^2} \, L_{t_0}(\beta) \leq \left.\frac{d}{dt}\right|_{t=t_0}
L_t(\beta) \leq \sqrt{a^2+b^2}\, L_{t_0}(\beta).$$
\end{proof}

Combining this with the inequality of Theorem \ref{thm:cusp-1},
we obtain the following theorem.

\begin{theorem}
  Let $L_t$ be the normalized length of a slope on a cusp torus
  $\partial U_t$ under a cone deformation, with derivative $\dot{L} =
  \frac{d}{dt} L_t$.  Then the change in $L_t$ is bounded by the
  following inequality:
  $$
  -\left(\sqrt{\frac{\sum b_i(t)}{2\mbox{Area}(\partial
        U_t)}}\right) \, L_t
	\;\leq\; \dot{L} \;\leq\;
  \left(\sqrt{\frac{\sum b_i(t)}{2\mbox{Area}(\partial U_t)}}\right)\,
  L_t.
  $$
\label{thm:norm-length-1}
\end{theorem}

\begin{proof}
  Solving for $\|\gamma'(t)\|$ in the inequality of Theorem
  \ref{thm:cusp-1}, we find
  $$\|\gamma'(t)\| \leq \sqrt{\frac{\sum b_j(t)}{2\mbox{Area}(\partial
      U_t)}}.$$
  The result is given by substituting this into equation
  (\ref{eqn:norm-length-L}).
\end{proof}

\subsubsection{Normalized height}

For our applications, the curve on the torus we are particularly
interested in is the one running orthogonal to the meridian.  When we
are dealing with link complements in $S^3$, as in a future section of
this paper, there is a well defined meridian on a cusp torus.
Otherwise, choose the shortest nontrivial simple closed curve on the
initial cusp torus to be the meridian (for our applications to knots
and links in $S^3$, these choices will agree).

We will be interested in estimating the length of the next shortest
nontrivial simple closed curve.  This length will be at least as long
as the length of the arc orthogonal to the meridian.  

At time $t$, let $p_t$ denote the geodesic arc perpendicular to the
meridian, running from meridian to meridian on $\partial U_t$.  This
arc will generally not be a closed curve.  Let $h_t$ be the normalized
length of $p_t$, that is, $h_t =
\rm{Length}(p_t)/\sqrt{\rm{Area}(\partial U_t)}$.  We refer to $h_t$
as the \emph{normalized height} of the torus at time $t$.  Denote its
derivative at time $t$ by $\dot{h}$.

\begin{lemma}
The normalized height $h_t$ satisfies
$$ -\sqrt{a^2+b^2}\,h_t \leq \dot{h} \leq
\sqrt{a^2+b^2}\,h_t.$$
\label{lemma:hdot-bound}
\end{lemma}

\begin{proof}
Let $\mu$ denote the meridian slope, with $L_t(\mu)$ its normalized
length at time $t$.  Then $h_t L_t(\mu) = 1$.  Thus
$$\frac{\dot{h}}{h_t} + \frac{\dot{L}(\mu)}{L_t(\mu)} = 0.$$
By Lemma \ref{lemma:ldot-bound} and Lemma \ref{lemma:a+b-gamma},
$$-\sqrt{a^2+b^2} \leq \frac{\dot{L}(\mu)}{L_t(\mu)} \leq
\sqrt{a^2+b^2}.$$ Hence
$ -\sqrt{a^2+b^2}\,h_t \leq \dot{h} \leq
\sqrt{a^2+b^2}\,h_t.$

\end{proof}

Again Lemma \ref{lemma:hdot-bound}, combined with the inequality of
Theorem \ref{thm:cusp-1}, gives
\begin{equation}
  -\left(\sqrt{\frac{\sum b_i}{2\mbox{Area}(\partial U_t)}}\right) \,
h_t \leq \dot{h} \leq
\left(\sqrt{\frac{\sum b_i}{2\mbox{Area}(\partial U_t)}}\right) \,
h_t
\label{eqn:norm-height}
\end{equation}

We simplify equation (\ref{eqn:norm-height}) by noting the area of
$U_t$ is given by the length of the meridian times the actual length
of the arc $p_t$.  Recall that the length of $p_t$ is
$\sqrt{\mbox{Area}(\partial U_t)}\,L_t(p_t) =
\sqrt{\mbox{Area}(\partial U_t)}h_t.$ So denoting the length of the
meridian of $\partial U_t$ by $m_t$, we have
$$\sqrt{\mbox{Area}(\partial U_t)} = m_t h_t.$$

Putting this into equation (\ref{eqn:norm-height}), we obtain:
\begin{theorem}
  Let $M$ be a hyperbolic 3--manifold with a cusp, and let $X$ be a
  hyperbolic 3--manifold which can be joined to $M$ by a smooth family
  of hyperbolic cone manifolds.  Let $h_t$ denote the normalized
  height at time $t$ of the cusp which deforms to the cusp of $M$, and
  let $m_t$ denote the length of its meridian.  Finally, denote by
  $\dot{h}$ the derivative of $h_t$.  Then
\begin{equation}
  -\frac{\sqrt{\sum b_i(t)}}{m_t\sqrt{2}} \leq \dot{h} \leq
\frac{\sqrt{\sum b_i(t)}}{m_t\sqrt{2}}
\label{eqn:norm-height-mer}
\end{equation}
\label{thm:norm-height-mer}
\end{theorem}

%% file: bdrybound.tex
In this section, we finish the proofs of Theorems \ref{thm:main-cusp}
and \ref{thm:main-norm-length} by
finding a parameterization of the cone deformation for which we may
bound the sum $\sum b_j$.

In order to bound $\sum b_j$, we need more explicit formulas for the
$b_j$.  Recall that $b_j$ was defined in Definition \ref{def:bj} as:
$$b_j = b_{V_j}(\eta_0, \eta_0).$$
Here $\eta_0$ is the real part of the 1--form $\omega_0$ in the $j$-th
component $V_j$ of a tubular neighborhood of the singular locus.  Thus
we begin by revisiting our decomposition of $\omega_0$ in $V_j$.

	
\subsection{Convex combinations of deformations}
	\label{subsec:convex}

Recall that corresponding to an infinitesimal deformation is the
1--form $\omega_0$.  We want $\omega_0$ to have certain desired
properties.

By equation (\ref{eqn:omega_0-decomp}), we can decompose $\omega_0 =
(\omega_0)_j$ in $V_j$ as
$$(\omega_0)_j =u_j\, \omega_m + (x_j+iy_j)\,\omega_\ell,$$
with $u_j$, $x_j$, and $y_j$ in $\RR$, and where $\omega_m$ and
$\omega_\ell$ were computed in \cite{hk:cone-rigid}.

In \cite{hk:univ}, it was shown that if the singular locus has just
one component, then the entire cone deformation can be parameterized
by the square of the cone angle $t=\alpha^2$.  In this case, $u=u_1$
was determined explicitly in \cite{hk:univ}: $u = -1/(4\alpha^2),$
where $\alpha$ is the cone angle at the core of the singular solid
torus $V=V_1$.

When we have multiple components of the singular locus, we may not
necessarily be able to parameterize the cone deformation in this
manner.  However, we do have the following information.

First, for each $j$ in $\{1, \dots, n\}$, locally we have a
deformation given by changing the $j$-th cone angle only and leaving
the others fixed.  This is a deformation with only one component of
the singular locus, so the methods of \cite{hk:univ} apply, and this
deformation can be parameterized by $t=\alpha_j^2$.  Then
\begin{equation*}
(\omega_0)_j = -\frac{1}{4\alpha_j^2}\, \omega_m +
(x_{jj}+iy_{jj})\omega_\ell.
\end{equation*}

Changing the angle in the $j$-th tube $V_j$ and leaving the others
fixed also affects the $k$-th tube $V_k$.  Since there is no change of
cone angle in this tube, $(\omega_0)_k$ can be expressed as
\begin{equation*}
  (\omega_0)_k = (x_{jk}+iy_{jk})\omega_\ell.
\end{equation*}

Any of these $j$ deformations may be scaled by a factor $s_j\geq 0$,
and we may take non-negative linear combinations.  These correspond to
new local deformations, in which the rates of change of cone angles
vary according to the choice of the $s_j$.  In a neighborhood of the
$j$-th component, write:
\begin{equation}
(\omega_0)_j = -\frac{s_j}{4\alpha_j^2}\,\omega_m + \sum_{k=1}^n (s_k
x_{kj} + i\,s_k y_{kj})\,\omega_\ell.
\label{eqn:omega_0-S}
\end{equation}

Given a point $\bar{s} = (s_1, s_2, \dots, s_n)$, we may use these
equations for $\omega_0$ to compute the forms $\eta_0$ explicitly, and
thus compute the boundary terms $b_j = b_j(\bar{s})$ explicitly.
Hodgson and Kerckhoff did the calculations without the extra $\bar{s}$
(page 382 of \cite{hk:univ}).  When the $\bar{s}$ is put in, we
obtain:
\begin{equation}
\frac{b_j(\bar{s})}{\mbox{Area}(\partial V_j)}= s_j^2\,c +
\left( (\sum_k s_k x_{kj})^2 + 
   (\sum_k s_k y_{kj})^2 \right)\,a +
  s_j(\sum_k s_k x_{kj}) \,b.
\label{eqn:b_j}
\end{equation}
Here the terms $a$, $b$, and $c$ are constants in terms of $R$ which
agree with those on page 383 of \cite{hk:univ}.  However, we are not
immediately concerned with their values.  More important, $a$ is
always negative, $b$ and $c$ always positive.

From the formula (\ref{eqn:b_j}) we obtain the following information.
\begin{enumerate}
  \item $b_j(\bar{s})$ is quadratic in $\bar{s}$.
  \item If any $s_j=0$ for some $\bar{s}$, then $b_j(\bar{s})\leq 0$.
\end{enumerate}

We also know that the sum of all the $b_j$'s is strictly positive,
provided $\bar{s} \neq (0, 0, \dots, 0)$, by equation
(\ref{eqn:b_j>0}), as well as the remark right after that equation.
We will use this information in \S\ref{subsec:Vnonempty}.


\subsection{Selecting local deformations}
  \label{subsec:Vnonempty}

Each choice of $\bar{s} = (s_1, \dots, s_n)$ corresponds to a local
deformation with 1--form $(\omega_0)_j$ expressed in $V_j$ as in
equation (\ref{eqn:omega_0-S}).
In terms of the deformation, varying $\bar{s}$ varies the rates at
which the cone angles are changing (instantaneously).  When
$\bar{s}=(1, \dots, 1)$, all cone angles are changing at the same
rate.  It might be convenient to let $\bar{s}=(1, \dots, 1)$, to
simplify calculations.  However, such a choice may not give us the
bounds on $\sum b_j$ that we need.

In particular, we would like to select $\bar{s}$ such that we may use
certain inequalities from \cite{hk:univ}.  (These inequalities appear
in the proof of Lemma \ref{lemma:b_j-bound} in this paper.)  In order
to use these inequalities, we will need to find $\bar{s}$ such that
each $b_j(\bar{s}) \geq 0$.

The existence of such an $\bar{s}$ has been discovered by
Hodgson and Kerckhoff, but their result is unpublished, so we include
it as Lemma \ref{lemma:Vnonempty} here.

\begin{lemma}[Hodgson--Kerckhoff]  There is an $\bar{s}\neq 0$ for which
$b_j(\bar{s})\geq 0$ for all $j=1, 2, \dots, n.$
\label{lemma:Vnonempty}
\end{lemma}

\begin{proof}
It suffices, by rescaling, to restrict $\bar{s}$ to the simplex
$$T=\{\bar{s}=(s_1, \dots,s_n) | \sum_j s_j = n, s_j \geq 0\}.$$

We first set up some notation.

Let the $j$-th vertex of $T$ be denoted $t_j$.  So $t_j = (0, \dots,
0, n, 0, \dots, 0)$ where the single non-zero value is in the $j$-th
location.  Notice that at $t_j$, $b_i(t_j)\leq 0$ for all $i\neq j$.
Because the sum of the $b_i$ is always strictly positive, we must have
$b_j(t_j)>0$.

Let the face of $T$ opposite $t_j$ be denoted $f_j$.  That is, on
$f_j$ the $j$-th coordinate is always $0$.  Notice that on $f_j$, $b_j
\leq 0$.  We will assume that $b_j < 0$ on $f_j$.  If not, consider
instead level sets of $b_j$ for values $\epsilon > 0$ arbitrarily
close to $0$.  Below, we will be finding points where all the $b_j$
are positive.  Taking a limit will at least give points where the
$b_j$ are nonnegative.

We will also assume level sets $b_j=0$ intersect transversely.  If
not, take an arbitrarily small perturbation of the quadratic functions
(keeping them quadratic).

First we note that the set where $b_j>0$ is an open neighborhood of
$t_j$ in $T$.  This is because on any line from $t_j$ to $f_j$, $b_j$
goes from positive to negative.  Because $b_j(s)$ is quadratic in $s$,
$b_j=0$ on this line in exactly one value.  Thus the level set $b_j=0$
separates $t_j$ from $f_j$.  It has a well--defined positive and
negative side.

Now, we prove the following sublemma by induction.

\begin{lemma}
  Consider the $k-1$ simplex $S_k$ of $T$ where $s_i=0$ for all $i>k$.
  There exists a nonempty open set $B_k \subseteq S_k$ where $b_j>0$
  for all $j\leq k$.  Furthermore, the set $E_k\subseteq S_k$ where
  $b_j=0$ for all $j<k$ is an odd number of points, each on the
  boundary of $B_k$.  
\label{lemma:sub-vnonempty}
\end{lemma}

Notice that in the above statement, when $k=n$, the simplex $S_n$ is
all of $T$.  Then if the statement is true, $B_n$ is an open set on
which all the $b_j$ are positive, so this will prove Lemma
\ref{lemma:Vnonempty}. 

An additional fact we will pick up from the proof is that the set
$E_n$ where $b_j=0$ for $j=1, 2, \dots, n-1$ is an odd number of
points.  (We will use this additional fact in the proof of Lemma
\ref{lemma:choose-S} below.)

\begin{proof}
We first prove Lemma \ref{lemma:sub-vnonempty} when $k=2$.  In this
case, $S_2$ is the 1--simplex from $t_1$ to $t_2$.  At $t_1$, $b_1>0$.
At $t_2$, $b_1 <0$.  Since $b_1$ is quadratic, somewhere on $S_2$ is a
single point $e_2$ where $b_1=0$.  Notice that $b_2(e_2)$ must be
positive.  This is because each $b_j(e_2)\leq 0$ for $j>2$ (since the
corresponding coordinates $s_j$ are all zero).  Since $b_1(e_2)$ is
zero, and the sum $\sum b_i >0$, this forces $b_2>0$.  Then in this
case, $E_2$ contains the single point $e_2$.  $B_2$ is nonempty, with
$e_2$ on its boundary.

\begin{figure}
  \input{figures/findv.pstex_t}
  \caption{Proof of Lemma \ref{lemma:Vnonempty} when $k=3$.}
  \label{fig:findv}
\end{figure}
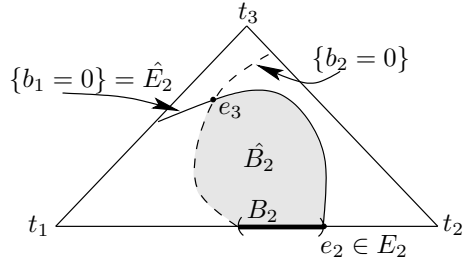

We will now prove the lemma by induction for $k+1$.  For reference, an
example of the case when $k+1=3$ is illustrated in Figure
\ref{fig:findv}.

$S_{k+1}$ is the $k$ subsimplex of $T$ with vertices $t_1, t_2, \dots,
t_{k+1}$, on which $s_{k+2}= \dots = s_n=0$.  It has as one boundary
face the $k-1$ simplex $S_k$, on which $s_{k+1}= s_{k+2}= \dots
=s_n=0$.  On any other face of $S_{k+1}$, $s_{k+2}=\dots=s_n=0$, and
also $s_j=0$ for some $j<k+1$.

By induction, we are assuming that on the face $S_k$ the set $B_k$,
where $b_j>0$ for all $j\leq k$ is a nonempty open set.  By
continuity, the set $\hat{B}_k$ in $S_{k+1}$ where the same $b_j>0$
must also be nonempty.  
The set $\hat{E}_k \subset S_{k+1}$ where $b_j=0$ for $j<k$ is a
$1$--manifold with boundary.  A neighborhood of $E_k$ in $\hat{E}_k$
lies in the boundary of $\hat{B}_k$.

Consider $\hat{E}_k$.  As a 1--manifold, it consists of closed
components and arcs.  The part of its boundary contained in $S_k$ is
the set $E_k$, so consists of an odd number of points.  Hence there
must be an odd number of arc components of $\hat{E}_k$ with exactly
one endpoint in $S_k$.

Consider the other boundary components of $\hat{E_k}$.  They must lie
on faces of $S_{k+1}$.  Recall the set $\{b_j=0\}$ does not intersect
any face $f_j$ where $s_j=0$.  Thus the only possible additional face
of intersection of $\hat{E_k}$ and $S_{k+1}$ is the face on which
$s_k=s_{k+2}=\dots=s_n=0$, spanned by vertices $t_1, \dots, t_{k-1}$
and $t_{k+1}$ (note this face is in $f_k$).  By induction (using the
induction hypothesis with vertices renumbered), there are an odd
number of intersection points of $\hat{E_k}$ on this face of
$S_{k+1}$.

Thus there must be an odd number of arcs of $\hat{E_k}$ running from
$S_k$ to the face $f_k$.  Along any of these arcs, $b_k$ goes from
positive (at a point of $E_k$ on $S_k$) to negative (since $b_k<0$ on
$f_k$).  Thus it must be zero at an odd number of interior points.
Each of these points is in $E_{k+1}$.  We may also have points of
$E_{k+1}$ coming from the intersection of $\{b_k=0\}$ with closed
components of $\hat{E}_k$, or with arc components of $\hat{E}_k$ with
both endpoints of the arc on the same face (either $S_k$ or $f_k$).
In the former case, $\{b_k=0\}$ must intersect any closed component an
even number of times, and in the latter case $b_k$ has the same sign
on both endpoints, so again $\{b_k=0\}$ intersects the arc an even
number of times.  Thus the total number of points in $E_{k+1}$ is odd.

At each point in $E_{k+1}$, we have $b_1=b_2=\dots=b_k=0$, and we have
$b_{k+2}$ through $b_n$ all less than or equal to zero.  Since the sum
of all the $b_j$ is positive, $b_{k+1}$ must be positive.  This
implies $B_{k+1}$ is a nonempty open set, and the points of $E_{k+1}$
are on its boundary.
\end{proof}

This completes the proof of Lemma \ref{lemma:Vnonempty}.
\end{proof}

Now, as mentioned above, the choice of $\bar{s}$ is tied to the rates
at which the cone angles are varying.  Lemma \ref{lemma:Vnonempty}
implies that a desired $\bar{s}$ exists, but gives us no further
information about $\bar{s}$.  In particular, some cone angle may reach
$2\pi$ before the others.  When this happens, the hyperbolic structure
in the corresponding tube about the singular locus is nonsingular.  We
may thus view the deformation from this point on as a deformation on a
cone manifold with one fewer component of the singular locus.  We do
this, obtaining deformations with fewer and fewer singular components,
until each cone angle has reached $2\pi$.


\subsection{Boundary bound}
\label{subsec:bd}

In this subsection, we will show that
there exists a parameterization of the cone deformation such that the
sum of the terms $b_j$ is of order $n$.  We first need two lemmas.
Lemma \ref{lemma:b_j-bound} gives us an initial bound on each $b_j$.
Lemma \ref{lemma:choose-S} will allow us to improve this bound.  It is
in the proof of Lemma \ref{lemma:b_j-bound}, giving the initial bounds
on the $b_j$, that we use the positivity result of Lemma
\ref{lemma:Vnonempty}.

Finally, in Theorem \ref{thm:bound-bdry-time} at the end of this
subsection, we put these lemmas together to record the main result of
the section: the bound on the sum $\sum b_j$.

\begin{lemma}
  Suppose we have a cone deformation given as a linear combination of
  deformations, as above, where cone angles go from $0$ to $2\pi$
  along the singular locus.  Suppose also that $S=(s_1, \dots, s_n)$
  is chosen such that for each time $t$, $b_j(S)\geq 0$ for all $j$.
  Suppose the tube radius $R$ is bounded below by $R_1 \geq 0.56$ for
  each $t$.  Then
  $$b_j(S) \leq s_j^2 \,C(R_1),$$
  where $C(R_1)$ is a function of $R_1$
  alone which approaches $0$ as $R_1$ approaches infinity.
\label{lemma:b_j-bound}
\end{lemma}

\begin{proof}
  First, since we're considering a fixed $V_j$, we will drop the
  subscript $j$ in our notation throughout the proof (i.e. $V_j=V,$
  $(\eta_0)_j = \eta_0,$ $s_j=s,$ etc).
  
  Modify a calculation (p. 382 and 383 of \cite{hk:univ}) to include
  $s$, and equation (17) of that paper can be written:
$$\frac{b_j(S)}{\mbox{Area}(\partial V)}\; \leq \;\frac{s^2}{8m^4},$$
where
$m = \alpha\sinh^2 R$ is the length of a meridian on $\partial U$.

We will bound the term $\mbox{Area}(\partial V)\, s^2/(8m^4)$.

Let $A=\mbox{Area}(\partial V)$.  Let $h$ denote the height of
$\partial V$, i.e. the length of an arc perpendicular to the
meridian.  So $h=\ell\cosh R$, where $\ell$ is the length of the
singular locus of the solid torus $V$, and $A=mh$.

Consider
$$\Phi = \frac{8m^4}{A} = \frac{8m^3}{h} =
\frac{8\alpha^3\sinh^3R}{\ell\cosh R}.$$
We will bound $\Phi$ below, therefore bounding its reciprocal above.

Let
$$\Upsilon(R) = \frac{3.3957 \tanh R }{2\cosh 2R}.$$
It was shown in
\cite{hk:univ} (see Theorem 4.4, Corollary 5.1, and the comments on
the multicusp case on p. 411), that $\alpha\ell \geq \Upsilon(R)$.

Then
$$\Phi = \frac{8\alpha^3\sinh^3R}{\ell\cosh R} =
\frac{4\alpha^3}{\ell}\,\frac{1}{\Upsilon(R)}\,g(R)$$
where
$$g(R) = \frac{2\Upsilon(R)\sinh^3R}{\cosh R} = 2\Upsilon(R)\tanh R\sinh^2
R = \frac{3.3957 \tanh^4 R}{1+\tanh^2 R}.$$

So $$\Phi \geq \frac{4\alpha^3}{\ell}\,\frac{1}{\alpha\ell}\,g(R) =
4\left(\frac{\alpha}{\ell}\right)^2g(R).$$

Note that $g(R)$ is increasing with $R$.  Hence since $R\geq R_1$,
$g(R)$ is bounded below by $g(R_1)$.  Thus we focus on the term
$\alpha/\ell$ in $\Phi$.

Define $u=\alpha/\ell$.  In our hyperbolic cone deformation, the
initial cone angle $\alpha$ is zero.  As $\alpha$ approaches $0$, $u$
approaches $L_0^2$, where $L_0$ is the normalized length of the
curve which becomes the meridian (i.e. bounds a singular solid disc)
under the cone deformation (see equation (37) of \cite{hk:univ}).  We
will estimate $u=(u_0 + \Delta u) = (L_0^2+\Delta u)$.

Reworking the calculations of \cite{hk:univ} to include the speed $s$,
and taking derivatives with respect to cone angle $\alpha$ rather than
time, we may rewrite Proposition 5.6 of that paper as:
$$\frac{1}{\alpha}\,\frac{du}{d\alpha} \geq -\frac{2(1+\tanh^2 R)}
{3.3957\tanh^3R} = -\frac{2\tanh R}{g(R)},$$
provided that $R\geq
R_1\geq 0.531$.

It is at this step, to use Proposition 5.6 of \cite{hk:univ}, that we
have required $S$ to be chosen so that $b_j(S) > 0$.

Thus
$$\Delta u \geq \int_{\alpha=0}^{\alpha=2\pi} -2\alpha\frac{\tanh
  R}{g(R)} d\alpha \geq
-\frac{\tanh R_1}{g(R_1)}\int_{0}^{2\pi} 2\alpha d\alpha =
-(2\pi)^2\frac{\tanh R_1}{g(R_1)}.$$

So
$$u^2 \geq (L_0^2 - (2\pi)^2\tanh(R_1)/g(R_1))^2,$$
provided the term
$L_0^2 - (2\pi)^2\tanh(R_1)/g(R_1)$ is non--negative.  We can ensure
this quantity is non--negative by ensuring $L_0^2$ is large.  This can
be done following \cite{hk:univ}.  The following result is essentially
equation (47) of that paper.

\begin{equation}
L_0^2 \geq \frac{2(2\pi)^2}{3.3957(1-\tanh R)}
  \exp\left( \int_1^{\tanh R} F(w) dw\right) = I(R),
  \label{eq:Lhat_int}
\end{equation}
provided the tube radius is at least $R$ throughout the deformation.
Here $F(w) = -(1+4w+6w^2+w^4)/((w+1)(1+w^2)^2)$ is an integrable
function.

The right hand side $I(R)$ of (\ref{eq:Lhat_int}) is increasing with
$R$, so it can be bounded below by a constant in terms of $R_1$.
$$L_0^2 \geq I(R_1)$$

Now if $R_1\geq 0.56$, then 
$I(R_1) - 2(2\pi)^2\tanh(R_1)/g(R_1) \geq 0$.

Putting this together,
$$\Phi \geq 4\left(I(R_1)-(2\pi)^2\frac{\tanh R_1}{g(R_1)}\right)^2
g(R_1).$$
The right hand side of the above equation is a constant,
$\tild{C}_0(R_1)$.

Thus we have
$$b_j(S) \leq \frac{A}{8m^4}\,s^2 =
  \frac{1}{\Phi}\,s^2 \leq \frac{s^2}{\tild{C}_0(R_1)} = s^2\,C(R_1).$$
\end{proof}


The result of Lemma \ref{lemma:b_j-bound} gives us a bound on the sum
of the $b_j$ in terms of $\bar{s}$.  If we could choose $\bar{s} = (1,
\dots, 1)$ (all cone angles increasing at the same rate), then Lemma
\ref{lemma:b_j-bound} would give a bound on $\sum b_j$ of order $n$.
However, our choice of $\bar{s}$ comes from Lemma
\ref{lemma:Vnonempty}.  By the proof of that lemma, we may assume only
that $\bar{s}$ is in the simplex $T=\{\bar{s}\, |\, \sum s_j=n,
s_j\geq 0\}$.  A priori, $\bar{s}$ could be close to a vertex of $T$,
say $(n, 0, \dots, 0)$, in which case Lemma \ref{lemma:b_j-bound}
gives only a bound of order $n^2$ on $\sum b_j$.

For our applications, we need the bound of order $n$.  We will obtain
this better estimate by again revisiting the choice of the point
$\bar{s}=(s_1, s_2, \dots, s_n)$.  This is done in Lemma
\ref{lemma:choose-S}.

\begin{lemma}
  There exists a point $S=(s_1, \dots, s_n)$ in $T$ such that each
  $b_j(S)\geq 0$, $j=1, \dots, n$, and for any index $j$ with $s_j
  \neq \min\{s_1, \dots, s_n\}$, we have $b_j(S)=0$.
\label{lemma:choose-S}
\end{lemma}

\begin{proof}

First we set up some notation.

\begin{notation}
  Let $B = \{\bar{s}\, |\, b_j(\bar{s}) \geq 0, j=1, \dots, n\}$ (so
  the interior of $B$ is the set $B_n$ of Lemma
  \ref{lemma:sub-vnonempty}).  Again we may normalize such that the
  points $\bar{s}$ lie in the simplex $T= \{(s_1,\dots,s_n)| \sum s_i
  = n, s_i\geq 0\}$.
  
  We define some sets in $B$.  For each $j=1, 2, \dots, n$, let $e_j$
  be the set of all points in $B$ such that $b_j>0$ and for $i\neq j$,
  $b_i=0$.  Note for $j=n$, $e_n$ corresponds to the set $E_n$ in
  Lemma \ref{lemma:sub-vnonempty}.  (So in particular, $e_n$ consists
  of an odd number of points.)
  
  Let $b_{\{i,j\}}$ be the set $b_{\{i,j\}} = B\cap (\cap_{k\neq i,j}
  \{b_k=0\})$.  Note $b_{\{i,j\}}$ has dimension $1$, with components with
  boundary in $e_i$ and $e_j$.
  
  In general, let $I$ be a $k$-element subset $I= \{i_1, i_2, \dots,
  i_k\} \subseteq \{1, \dots, n\}$.  Let $b_I = B\cap(\cap_{j\notin I}
  \{b_j=0\})$.  Notice $b_I$ is bounded by components $b_{I_j}$ where
  $j\in I$, and $I_j$ is the $(k-1)$-element set obtained by removing
  $j$ from $I$.
  
  Now we define some sets in the simplex $T = \{\bar{s}| \sum
  s_i=n\}$.  Let $T_i$ be the region in which $s_i\leq s_j$ for all
  indices $j$.  $T_i$ has dimension $n-1$.  Let $T_{\{i,j\}} = T_i
  \cap T_j$.  Thus each $T_i$ is bounded by faces $T_{\{i,j\}}$ where
  $j$ ranges over all elements of $\{1, \dots, n\}-\{i\}$, as well as
  faces on $\partial T$.  Generally,
  let $I$ be a $k$-element subset of $\{1, \dots, n\}$.  Let $T_I =
  \cap_{j\in I} T_j$.  Note the boundary faces of $T_I$ are the faces
  given by $T_{I\cup\{k\}}$, where $k$ ranges over all elements of
  $\{1, \dots, n\} - I$, as well as faces in $\partial T$.  See Figure
  \ref{fig:lemma-fix} for an illustration when $n=3$.
\end{notation}

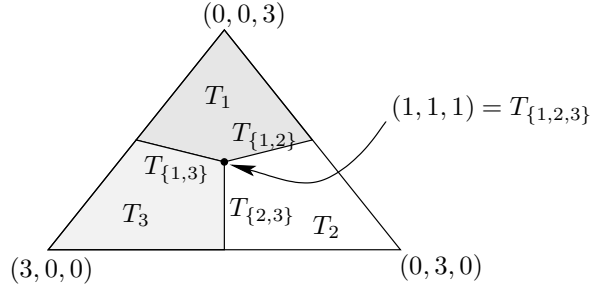
\begin{figure}
	\input{figures/fix-lemma.pstex_t}
	\caption{Faces $T_I$ illustrated when $n=3$.}
	\label{fig:lemma-fix}
\end{figure}

We will make transversality assumptions that allow us to ignore the
boundary faces of $T_I$ on $\partial T$.  In particular, assume that
$B$ lies in the interior of $T$, so the sets $b_I$ do not meet
$\partial T$.  If this is not the case, we may apply the argument
below to level sets $b_I = \epsilon$, then take the limit as $\epsilon
\to 0$.  Since $B$ is closed, we will obtain the desired result.
Similarly, we will assume that the $b_I$'s and the $T_I$'s all
intersect transversely.

We will use the transversality assumption that $B$ does not meet
$\partial T$ in the arguments below.  In particular, we will be
considering the intersection of subsets of $B$ with boundary faces of
$T_I$.  Since $B$ does not meet $\partial T$, this allows us to ignore
the faces of $T_I$ which lie on $\partial T$ in our argument.


Now, note the lemma will be immediately true if we can find some index
$j$ such that $e_j$ is in the region $T_j$: for on $e_j$, $b_i=0$ for
all $i\neq j$, and in $T_j$, $s_j$ is minimum.  So assume this never
happens.

Generally, if for some $k$-element set $I \subseteq \{1, \dots, n\}$,
$b_I \cap T_I$ is nonempty, we obtain a point which satisfies the
lemma.  So we will assume this also never happens.  We will show that
in this case, the point $(1, 1, \dots, 1)$ is in $B$, thus concluding
the proof.  To prove this, we need another sublemma.

\begin{lemma}
Let $\#(\cdot)$ denote the number of elements of a finite set.
Assuming that 
each 
$b_I \cap T_I$ is empty, 
$$\sum_I \#(b_{I\cup\{n\}} \cap \partial T_I) \equiv \sum_{J}
\#(b_{J\cup\{n\}} \cap \partial T_J) \equiv 1 \mod 2,$$
where the first
sum is over all $(k+1)$-element subsets $I$ of $\{1, \dots, n-1\}$,
and the second sum is over all $k$-element subsets $J$ of $\{1, \dots,
n-1\}$.
\label{lemma:inductive-S}
\end{lemma}

Check that the statement of the lemma makes sense:  Note if $J$ is a
$k$-element set, then $b_{J\cup\{n\}}$ is $k$-dimensional.  $T_J$ is
$(n-k)$ dimensional, so its boundary faces are $(n-k-1)$-dimensional.
Thus the intersections in the statement of the lemma are
$0$-dimensional manifolds, consisting of a finite number of
points.

Also notice that when $k=(n-2)$ in the above lemma, there exists only
one $(k+1)$-element subset $I$, namely $I=\{1, \dots, n-1\}$.  So
$b_{I\cup\{n\}}$ is the whole set $B$.  On the other hand, the
boundary $\partial T_I$ in this case must be the point in the
intersection of all the $T_i$, or the point $(1, \dots, 1)$.  Thus the
first sum counts the number of times the point $(1, \dots, 1)$ appears
in $B$.  Hence showing this number is odd will finish the
proof of Lemma \ref{lemma:choose-S}.

\begin{proof}
We prove Lemma \ref{lemma:inductive-S} by induction on $k$.

We start by showing what essentially is the case $k=0$.  That is, we
show
$$\sum_{i=1}^{n-1} \#(b_{\{i,n\}} \cap \partial T_i) \equiv
\sum_{i=1}^{n-1} \#(e_n \cap T_i) \mod 2.$$

Since we are assuming $e_n \cap T_n = \emptyset$, the second sum in
the above equation is just the number of points of $e_n$ (since $e_n$
is disjoint from each $T_i \cap T_j$ by transversality).  That number
of points is odd by Lemma \ref{lemma:sub-vnonempty}, so this will give
us the basis step of our induction.

Consider $b_{\{i, n\}} \cap T_i$.  This is a 1--manifold, consisting
of closed curves and arcs.  Each arc has two endpoints.  So
$$\sum_{i=1}^{n-1} \#\partial(b_{\{i,n\}} \cap T_i) \equiv 0 \mod 2$$

On the other hand,
$$\sum_{i=1}^{n-1} \#\partial(b_{\{i,n\}} \cap T_i) =
\sum_{i=1}^{n-1} \#(\partial b_{\{i,n\}} \cap T_i) +
\sum_{i=1}^{n-1} \#(b_{\{i,n\}} \cap \partial T_i).$$
So $\sum_{i\neq n} \#(b_{\{i,n\}} \cap \partial T_i)$ has the same
parity as $\sum_{i\neq n} \#(\partial b_{\{i,n\}} \cap T_i)$.

Consider $\partial b_{\{i,n\}} \cap T_i$.  Note $\partial b_{\{i,n\}}$
consists of points in $e_i$ and $e_n$.  By assumption, $e_i$ is never
in $T_i$.  Thus
$$\sum_{i=1}^{n-1} \#(\partial b_{\{i,n\}} \cap T_i)
= \sum_{i=1}^{n-1} \#(e_n \cap T_i) = \#(e_n).$$
So $\sum_{i\neq n} \#(b_{\{i,n\}} \cap \partial T_i)$ is odd.

Our proof for general $k$ follows essentially the same lines as the
proof above.

Suppose $1\leq k \leq n-2$.  We need to show
$$\sum_I \#(b_{I\cup\{n\}} \cap \partial T_I) \equiv \sum_{J}
\#(b_{J\cup\{n\}} \cap \partial T_J) \mod 2,$$
where the first sum is over $(k+1)$-element subsets of $\{1, \dots,
n-1\}$ and the second is over $k$-element subsets.

Fix some $(k+1)$-element subset $I$, and again consider
$b_{I\cup\{n\}} \cap T_I$.  This is a 1--manifold.  Its components are
closed curves and arcs.  Again the number of endpoints of
these arcs is even.  So when we sum over all such $I$, 
$$\sum_I \#\partial(b_{I\cup\{n\}} \cap T_I) \equiv 0 \mod 2.$$

On the other hand, 
$$\sum_I \#\partial(b_{I\cup\{n\}} \cap T_I) =
\sum_I \#(\partial b_{I\cup\{n\}} \cap T_I) +
\sum_I \#(b_{I\cup\{n\}} \cap \partial T_I).$$

Thus we consider $\sum_I \#(\partial b_{I\cup\{n\}} \cap T_I).$ We
want it to have the same parity as $\sum_J \#(b_{J\cup\{n\}} \cap
\partial T_J)$, where $J$ ranges over all $k$-element subsets of $\{1,
\dots, n-1\}$.  We show now that the two sums are actually equal,
which concludes the proof.

Consider $\partial b_{I\cup \{n\}}$.  This consists of $b_I$, as well
as things of the form $b_{{I_j}\cup\{n\}}$ where $I_j$ is the
$k$-element set obtained from $I$ by removing the element $j\in I$.
Since by assumption, $b_I \cap T_I$ is empty, $\partial b_{I\cup\{n\}}
\cap T_I$ is the union of sets $b_{{I_j}\cup\{n\}} \cap T_I$, where
$j$ ranges over all elements of $I$.

Then
\begin{equation}
\sum_I \#(\partial b_{I\cup\{n\}} \cap T_I) = \sum_I \sum_{j\in I}
\#(b_{{I_j}\cup\{n\}} \cap T_I).
\label{eqn:count1}
\end{equation}

We will change our method of counting in (\ref{eqn:count1}) above.
Replace the terms $I_j$ with $k$-element subsets $J$ of $\{1, \dots,
n-1\}$.  We will run through all such $J$.  Then $I$ will run through
the $(k+1)$-element subsets given by $J\cup\{i\}$, where $i$ is some
element of $\{1, \dots, n-1\} - J$.  So the right hand side of
equation (\ref{eqn:count1}) is equal to
$$\sum_J \sum_{i\in\{1,\dots,n-1\}-J} \#(b_{J\cup\{n\}} \cap
T_{J\cup\{i\}}).$$

On the other hand, we know the boundary $\partial T_J$ is equal to the
union of all faces $T_{J\cup\{j\}}$, where $j$ ranges over all
elements of $\{1, \dots, n\}-J$.  The only such face we are missing in
our sum above is $T_{J\cup\{n\}}$.  But by assumption, the
intersection of $b_{J\cup\{n\}}$ with $T_{J\cup\{n\}}$ is trivial.
Thus
$$\sum_J \sum_{i\in\{1,\dots,n-1\}-J} \#(b_{J\cup\{n\}} \cap
T_{J\cup\{i\}}) = \sum_J \#(b_{J\cup\{n\}} \cap \partial T_J).$$
\end{proof}

This concludes the proof of Lemma \ref{lemma:choose-S}.
\end{proof}


Lemmas \ref{lemma:b_j-bound} and \ref{lemma:choose-S} together give
the desired bound on $\sum b_j$:

\begin{theorem}
  Let $X$ be a hyperbolic manifold with $n+k$ cusps connected to a
  manifold $M$ with $k$ cusps by a cone deformation.  (So in
  particular, there are $n$ components of the singular locus, and cone
  angles at each component go from $0$ to $2\pi$.)  Suppose during the
  deformation, the tube radius about the singular locus remains larger
  than $R_1 \geq 0.56$.
	Then the deformation can be parameterized so that at each time: $\min
	s_j =1$ and the sum of the boundary terms satisfies
  $$\sum_{j=1}^n b_j \leq n\,C(R_1).$$
  Moreover, we can assume the
  deformation has reached the hyperbolic structure on $M$ by time
  $t=(2\pi)^2$.
  \label{thm:bound-bdry-time}
\end{theorem}

\begin{proof}
  
  By Lemma \ref{lemma:choose-S}, we may find a value $S$ in $T$ such
  that each $b_j(S) \geq 0$ and such that if for any index $j$, $s_j
  \neq \min\{s_1, \dots, s_n\}$, then $b_j(S)=0$.  Denote the minimum
  value of the $s_j$ by $s$: $s=\min\{s_1,\dots, s_n\}$.  Then Lemma
  \ref{lemma:b_j-bound} implies that
\begin{equation}
  \sum_{j=1}^n b_j = \sum_{\{j|b_j\neq 0\}} b_j \leq \sum_{\{j|b_j\neq
  0\}} s^2\,C(R_1) \leq n s^2\,C(R_1).
\label{eqn:sum-b_j}
\end{equation}

First, we need to check that $s$ is not equal to zero.  Suppose not.
That is, suppose $s=0$.  Then by (\ref{eqn:sum-b_j}), $\sum_{j=0}^n
b_j = 0$.  Recall that by equation (\ref{eqn:b_j>0}) and the remark
right after it, the sum of all $b_j$'s equal zero only if the entire
deformation is trivial.  But this is possible only if each $s_j=0$.
Since the point $(0, 0, \dots, 0)$ is not in the simplex $T$, this is
impossible.

Now, recall the terms $s_j$ affect the rate of change of the
deformation.  Thus we may rescale them without changing the tube
radius.  For the point $S$ in $T$, divide each $s_j$ by $s$.
Replace $S$ by this rescaled point, and denote the rescaled
point by the same notation: $S=(s_1, \dots, s_n)$.  So each $s_j\geq
1$.  Also, $s=\min\{s_j\}=1$.

Thus equation (\ref{eqn:sum-b_j}) gives:
$$\sum_{j=1}^n b_j \leq n\,C(R_1).$$
  
Finally, to notice that the deformation is complete by time
$t=(2\pi)^2$ as claimed, recall that if $s_j=1$ throughout the
deformation, then the cone angle is changing at a rate at which
$\alpha_j$ will reach $2\pi$ at time $t=(2\pi)^2$.  (See the
discussion at the beginning of \S\ref{subsec:convex}.)  Since $s_j\geq
1$ always, the rate of change may only be faster.  Thus each cone
angle reaches $2\pi$ by the time $t=(2\pi)^2$.
\end{proof}

When we combine Theorem \ref{thm:bound-bdry-time} with the results of
Theorem \ref{thm:cusp-1} and Theorem \ref{thm:norm-length-1}, we have
completed the proofs of Theorems \ref{thm:main-cusp} and
\ref{thm:main-norm-length}.

%% file: figures/findv.pstex_t
\begin{picture}(0,0)%
\includegraphics{figures/findv.pstex}%
\end{picture}%
\setlength{\unitlength}{3947sp}%
\begingroup\makeatletter\ifx\SetFigFont\undefined%
\gdef\SetFigFont#1#2#3#4#5{%
  \reset@font\fontsize{#1}{#2pt}%
  \fontfamily{#3}\fontseries{#4}\fontshape{#5}%
  \selectfont}%
\fi\endgroup%
\begin{picture}(2755,1680)(892,-1191)
\put(2346,342){\makebox(0,0)[lb]{\smash{{\SetFigFont{10}{12.0}{\familydefault}{\mddefault}{\updefault}{\color[rgb]{0,0,0}$t_3$}%
}}}}
\put(3632,-990){\makebox(0,0)[lb]{\smash{{\SetFigFont{10}{12.0}{\familydefault}{\mddefault}{\updefault}{\color[rgb]{0,0,0}$t_2$}%
}}}}
\put(2221,-253){\makebox(0,0)[lb]{\smash{{\SetFigFont{10}{12.0}{\familydefault}{\mddefault}{\updefault}{\color[rgb]{0,0,0}$e_3$}%
}}}}
\put(1028,-994){\makebox(0,0)[lb]{\smash{{\SetFigFont{10}{12.0}{\familydefault}{\mddefault}{\updefault}{\color[rgb]{0,0,0}$t_1$}%
}}}}
\put(2399,-909){\makebox(0,0)[lb]{\smash{{\SetFigFont{10}{12.0}{\familydefault}{\mddefault}{\updefault}{\color[rgb]{0,0,0}$B_2$}%
}}}}
\put(2810, 50){\makebox(0,0)[lb]{\smash{{\SetFigFont{10}{12.0}{\familydefault}{\mddefault}{\updefault}{\color[rgb]{0,0,0}$\{b_2=0\}$}%
}}}}
\put(2378,-589){\makebox(0,0)[lb]{\smash{{\SetFigFont{10}{12.0}{\familydefault}{\mddefault}{\updefault}{\color[rgb]{0,0,0}$\hat{B_2}$}%
}}}}
\put(907,-83){\makebox(0,0)[lb]{\smash{{\SetFigFont{10}{12.0}{\familydefault}{\mddefault}{\updefault}{\color[rgb]{0,0,0}$\{b_1=0\} = \hat{E_2}$ }%
}}}}
\put(2861,-1127){\makebox(0,0)[lb]{\smash{{\SetFigFont{10}{12.0}{\familydefault}{\mddefault}{\updefault}{\color[rgb]{0,0,0}$e_2 \in E_2$}%
}}}}
\end{picture}%

%% file: figures/fix-lemma.pstex_t
\begin{picture}(0,0)%
\includegraphics{figures/fix-lemma.pstex}%
\end{picture}%
\setlength{\unitlength}{3947sp}%
\begingroup\makeatletter\ifx\SetFigFont\undefined%
\gdef\SetFigFont#1#2#3#4#5{%
  \reset@font\fontsize{#1}{#2pt}%
  \fontfamily{#3}\fontseries{#4}\fontshape{#5}%
  \selectfont}%
\fi\endgroup%
\begin{picture}(3046,1766)(362,-1246)
\put(1744,-822){\makebox(0,0)[lb]{\smash{{\SetFigFont{10}{12.0}{\familydefault}{\mddefault}{\updefault}{\color[rgb]{0,0,0}$T_{\{2,3\}}$}%
}}}}
\put(1761,-348){\makebox(0,0)[lb]{\smash{{\SetFigFont{10}{12.0}{\familydefault}{\mddefault}{\updefault}{\color[rgb]{0,0,0}$T_{\{1,2\}}$}%
}}}}
\put(2267,-933){\makebox(0,0)[lb]{\smash{{\SetFigFont{10}{12.0}{\familydefault}{\mddefault}{\updefault}{\color[rgb]{0,0,0}$T_2$}%
}}}}
\put(2758,-188){\makebox(0,0)[lb]{\smash{{\SetFigFont{10}{12.0}{\familydefault}{\mddefault}{\updefault}{\color[rgb]{0,0,0}$(1,1,1)= T_{\{1,2,3\}}$}%
}}}}
\put(1208,-558){\makebox(0,0)[lb]{\smash{{\SetFigFont{10}{12.0}{\familydefault}{\mddefault}{\updefault}{\color[rgb]{0,0,0}$T_{\{1,3\}}$}%
}}}}
\put(2809,-1183){\makebox(0,0)[lb]{\smash{{\SetFigFont{10}{12.0}{\familydefault}{\mddefault}{\updefault}{\color[rgb]{0,0,0}$(0,3,0)$}%
}}}}
\put(362,-1201){\makebox(0,0)[lb]{\smash{{\SetFigFont{10}{12.0}{\familydefault}{\mddefault}{\updefault}{\color[rgb]{0,0,0}$(3,0,0)$}%
}}}}
\put(1571,400){\makebox(0,0)[lb]{\smash{{\SetFigFont{10}{12.0}{\familydefault}{\mddefault}{\updefault}{\color[rgb]{0,0,0}$(0,0,3)$}%
}}}}
\put(1590,-116){\makebox(0,0)[lb]{\smash{{\SetFigFont{10}{12.0}{\familydefault}{\mddefault}{\updefault}{\color[rgb]{0,0,0}$T_1$}%
}}}}
\put(1072,-857){\makebox(0,0)[lb]{\smash{{\SetFigFont{10}{12.0}{\familydefault}{\mddefault}{\updefault}{\color[rgb]{0,0,0}$T_3$}%
}}}}
\end{picture}%

%% file: meridian.tex
We now show that the length of the shortest nontrivial simple closed
curve on the cusp torus is bounded throughout the cone deformation,
provided the tube radius about the singular locus stays large.  Recall
here that we are calling this shortest curve the meridian of the
torus.  The following lemma is the last piece needed to complete the
proof of Theorem \ref{thm:main-norm-height}.

\begin{lemma}
Suppose $M$ is a hyperbolic cone manifold with one cusp, and with
singular locus $\Sigma$.  Suppose the tube radius $R$ about the
singular locus is at least $\log(2/\sqrt{3})$.  Then there exists a
horoball neighborhood $U$ of the cusp such that the meridian of
$\partial U$ has length at least $1-e^{-2R}$.
  \label{lemma:meridian-bound}
\end{lemma}

\begin{proof}
Let $X = M-\Sigma$.  $X$ admits an incomplete hyperbolic metric.
Expand the tube $V = \cup V_i$ about the
singular locus $\Sigma$ until the tube hits itself.  Each $V_i$ has
the same tube radius $R$.  Next expand the horoball $U$ about
the cusp until it either becomes tangent to itself or to $V$.

We obtain a region $\Omega$ by continuing to expand $U$ as follows.
First, choose a particular lift $U_a$ of $U$ in the universal cover.
Because $U$ is embedded in $\HH^3$ by the developing map on the
universal cover $\widetilde{X}$ of $X$, we can assume the developing
map is one to one on $U_a$ and takes $\partial U_a$ to the plane $z=1$
in $\HH^3$.

Take $\Omega= U_a \cup S$, where $S$ is the set consisting of all
points $x$ in $\tild{X}$ such that the distance $\mbox{dist}(x,
\partial U_a) < \mbox{dist}(x, \partial U_b)$ for $U_b$ any other lift
of $U$, and also such that the distance $\mbox{dist}(\phi(x), \partial
U) < \mbox{dist}(\partial U,\Sigma)$, where $\phi$ is the covering
projection. Note that $\Omega$ is embedded in $\HH^3$ by the
developing map.  That is, the developing map is one to one on
$\Omega$.

Let $p$ be a point on $\partial U$ closest to $\Sigma$, and let $p_0$
and $p_1$ be lifts of $p$ on $\partial U_a$ such that the projection
of the path from $p_0$ to $p_1$ is shortest on $\partial U$.  Let
$a_0$ be the point in the universal cover on a lift of $\Sigma$ such
that the path from $p_0$ to $a_0$ has length equal to the minimum
distance from $p_0$ to $\Sigma$.  Similarly, let $a_1$ be the point on
a lift of $\Sigma$ such that the path from $p_1$ to $a_1$ has length
equal to the minimum distance from $p_1$ to $\Sigma$.  Note the
interiors of the paths between $p_0$ and $a_0$ and between $p_1$ and
$a_1$ are contained in $\Omega$. 

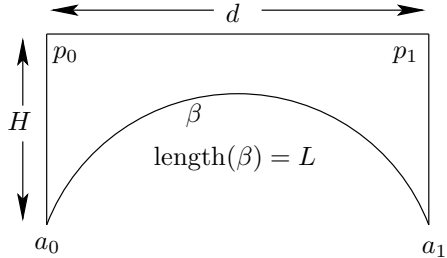
\begin{figure}
  \begin{center}
  \input{figures/area_quad.pstex_t}
  \end{center}
  \caption{The top edge of the quadrilateral maps to the shortest path
  on $\partial U$ under the covering projection. }
  \label{fig:area-quad}
\end{figure}

Then note the interior of the geodesic $\beta$ from $a_0$ to $a_1$ is
also in $\Omega$, along with all the points in the quadrilateral with
sides given by $\beta$ and the paths from $p_0$ to $p_1$, from $p_0$
to $a_0$, and from $p_1$ to $a_1$.  This is because, first, all points
$x$ in this quadrilateral (except the endpoints of $\beta$) have
$\mbox{dist}(\phi(x), \partial U) < \mbox{dist}(\partial U, \Sigma)$.
Also, if any $x$ were closer to some other lift $U_b$ of $U$, then
there would be points in the quadrilateral equidistant to two
different lifts of $U$.  All points equidistant to two different lifts
of $U$ lie on some totally geodesic planes of $\HH^3$.  Hence some
geodesic must intersect $\beta$ but not the paths from $p_0$ to $a_0$,
$p_1$ to $a_1$, or the path on the horosphere $\partial U_a$ from
$p_0$ to $p_1$.  This is impossible.

Now let $H$ be the length of the paths from $p_0$ to $a_0$ and from
$p_1$ to $a_1$.  Let $L$ be the length of $\beta$.  We know that
$H\geq R$ and $L\geq 2R$ since the distance between any two components
of the lift of $\Sigma$ is greater than or equal to $2R$.  We let $d$
be the length of the path on the horosphere $\partial U_a$ between
$p_0$ and $p_1$.  Since $\Omega$ is embedded by the developing map, we
can assume the quadrilateral is of the form shown in Figure
\ref{fig:area-quad}.

Then solving for $d$ in terms of $H$ and $L$, $d= 2e^{-H}\sinh(L/2)$.

Case 1:  $H=R$.  Then $d \geq 2e^{-R}\sinh(R) = 1-e^{-2R}$.

Case 2: $H>R$ Then $U$ must be tangent to itself, say at a point $q$,
and not necessarily tangent to the tube about $\Sigma$.  Then the
image of $\Omega$ under the developing map must contain some portion
of the region above the geodesic surface equidistant from a horizontal
plane (the lift of $\partial U_a$) and some other horosphere.

This geodesic surface is some portion of the (Euclidean) hemisphere
with radius $1$, tangent to the lift of $\partial U_a$ at a lift $q_0$
of $q$.  See Figure \ref{fig:d-long}.  The image of $\Omega$ under the
developing map contains a neighborhood of points on this hemisphere
about $q_0$.  If it were to contain all the points on the hemisphere
between the dotted lines in Figure \ref{fig:d-long} (that is, points
whose projection to $\CC$ have Euclidean distance at most $1/2$ from
the projection of $q_0$ to $\CC$), then $d\geq 1$.

Let $x$ be one of these points.  That is, the projection of $x$ to
$\CC$ has Euclidean distance at most $1/2$ from the projection of
$q_0$ to $\CC$.  Note $x$ will lie in $\Omega$ provided
$\mbox{dist}(\phi(x), \partial U) < \mbox{dist}(\partial U, \Sigma) =
H$, or provided $H$ is large enough.  But now, $\mbox{dist}(\phi(x),
\partial U) = \log( 1/\mbox{height}(x))$, where $\mbox{height}(x)$
denotes the $z$--coordinate of $x$ in $\HH^3$.  Thus this distance is
at most $\log(2/\sqrt{3})$.

\begin{figure}
\input{figures/d-long.pstex_t}	
\caption{}
\label{fig:d-long}
\end{figure}

Since $H>R$, we can ensure $d\geq 1$ by restricting to $R\geq
\log(2/\sqrt{3})$.
\end{proof}

When we plug the result of Lemma \ref{lemma:meridian-bound} into
Theorem \ref{thm:norm-height-mer}, along with the bound on $\sum b_j$
of the previous section, we may complete the proof of Theorem
\ref{thm:main-norm-height}.

\medskip

\noindent\textbf{Theorem \ref{thm:main-norm-height}.} \emph{
  Let $M=X_\tau$ be a hyperbolic 3--manifold with a cusp, and let
  $X=X_0$ be a hyperbolic 3--manifold which can be joined to $M$ by a
  smooth family of hyperbolic cone manifolds $X_t$ with tube radius at
  least $R_1 \geq 0.56$, and with $n$ components of the singular
  locus. Let $U_t$ be a horoball neighborhood about
  the cusp.  Let $h(M)$ denote the normalized height of the cusp torus
  $\partial U_\tau$.  Let $h(X)$ denote the normalized height of
  $\partial U_0$.  Then the
	hyperbolic cone deformation can be parameterized such that the
	change in normalized height is bounded in
  terms of $R_1$ alone:
$$-(2\pi)^2\frac{\sqrt{n\,C(R_1)}}{(1-e^{-2R_1})\sqrt{2}} \leq
 h(M) - h(X) \leq
(2\pi)^2\frac{\sqrt{n\,C(R_1)}}{(1-e^{-2R_1})\sqrt{2}}.$$
}

\begin{proof}
  Theorem \ref{thm:norm-height-mer}, combined with Theorem
  \ref{thm:bound-bdry-time} and Lemma \ref{lemma:meridian-bound}
  imply that
$$
  -\frac{\sqrt{n\,C(R_1)}}{(1-e^{-2R_1})\sqrt{2}} \leq \dot{h} \leq
\frac{\sqrt{n\,C(R_1)}}{(1-e^{-2R_1})\sqrt{2}}
$$
and that the total time of the deformation is $t=(2\pi)^2$.  We
integrate the above inequality over the deformation.  The bounds on
the left and right are independent of $t$.  Thus we obtain the
conclusion of the theorem.
\end{proof}

%% file: figures/area_quad.pstex_t
\begin{picture}(0,0)%
\includegraphics{figures/area_quad.pstex}%
\end{picture}%
\setlength{\unitlength}{3947sp}%
\begingroup\makeatletter\ifx\SetFigFont\undefined%
\gdef\SetFigFont#1#2#3#4#5{%
  \reset@font\fontsize{#1}{#2pt}%
  \fontfamily{#3}\fontseries{#4}\fontshape{#5}%
  \selectfont}%
\fi\endgroup%
\begin{picture}(2672,1637)(941,-1800)
\put(1126,-1711){\makebox(0,0)[lb]{\smash{{\SetFigFont{10}{12.0}{\familydefault}{\mddefault}{\updefault}{\color[rgb]{0,0,0}$a_0$}%
}}}}
\put(3556,-1736){\makebox(0,0)[lb]{\smash{{\SetFigFont{10}{12.0}{\familydefault}{\mddefault}{\updefault}{\color[rgb]{0,0,0}$a_1$}%
}}}}
\put(3376,-511){\makebox(0,0)[lb]{\smash{{\SetFigFont{10}{12.0}{\familydefault}{\mddefault}{\updefault}{\color[rgb]{0,0,0}$p_1$}%
}}}}
\put(2326,-286){\makebox(0,0)[lb]{\smash{{\SetFigFont{10}{12.0}{\familydefault}{\mddefault}{\updefault}{\color[rgb]{0,0,0}$d$}%
}}}}
\put(1241,-511){\makebox(0,0)[lb]{\smash{{\SetFigFont{10}{12.0}{\familydefault}{\mddefault}{\updefault}{\color[rgb]{0,0,0}$p_0$}%
}}}}
\put(956,-966){\makebox(0,0)[lb]{\smash{{\SetFigFont{10}{12.0}{\familydefault}{\mddefault}{\updefault}{\color[rgb]{0,0,0}$H$}%
}}}}
\put(2076,-916){\makebox(0,0)[lb]{\smash{{\SetFigFont{10}{12.0}{\familydefault}{\mddefault}{\updefault}{\color[rgb]{0,0,0}$\beta$}%
}}}}
\put(1876,-1186){\makebox(0,0)[lb]{\smash{{\SetFigFont{10}{12.0}{\familydefault}{\mddefault}{\updefault}{\color[rgb]{0,0,0}$\mbox{length}(\beta) = L$}%
}}}}
\end{picture}%

%% file: figures/d-long.pstex_t
\begin{picture}(0,0)%
\includegraphics{figures/d-long.pstex}%
\end{picture}%
\setlength{\unitlength}{3947sp}%
\begingroup\makeatletter\ifx\SetFigFont\undefined%
\gdef\SetFigFont#1#2#3#4#5{%
  \reset@font\fontsize{#1}{#2pt}%
  \fontfamily{#3}\fontseries{#4}\fontshape{#5}%
  \selectfont}%
\fi\endgroup%
\begin{picture}(1923,1348)(633,-856)
\put(2153,-573){\makebox(0,0)[lb]{\smash{{\SetFigFont{10}{12.0}{\familydefault}{\mddefault}{\updefault}{\color[rgb]{0,0,0}$\sqrt{3}/2$}%
}}}}
\put(1580,155){\makebox(0,0)[lb]{\smash{{\SetFigFont{10}{12.0}{\familydefault}{\mddefault}{\updefault}{\color[rgb]{0,0,0}$q_0$}%
}}}}
\put(1559,384){\makebox(0,0)[lb]{\smash{{\SetFigFont{10}{12.0}{\familydefault}{\mddefault}{\updefault}{\color[rgb]{0,0,0}$1$}%
}}}}
\put(1303,-52){\makebox(0,0)[lb]{\smash{{\SetFigFont{10}{12.0}{\familydefault}{\mddefault}{\updefault}{\color[rgb]{0,0,0}$x$}%
}}}}
\end{picture}%

%% file: hyp-knot.tex

We will now apply the results of the previous sections to determine
bounds on the shapes of cusps of knot complements.  In particular, let
$K$ be a hyperbolic knot with a prime, twist reduced diagram, where
prime and twist reduced are defined in Definitions \ref{def:prime} and
\ref{def:twred}.  We will
bound the normalized height of the cusp of $S^3-K$, using Theorem
\ref{thm:main-norm-height}.  We need the following:
\begin{enumerate}
  \item An initial hyperbolic manifold $X$ with $n+1$ cusps.
  \item A cone deformation deforming $X$ to the knot complement
    $S^3-K$.
  \item A lower bound on the tube radius ($R\geq R_1\geq 0.56$).
\end{enumerate}

Given these pieces, we will obtain bounds on
normalized height.
  
\subsection{An initial manifold}
\label{subsec:init}

\begin{figure}
\begin{center}
\includegraphics{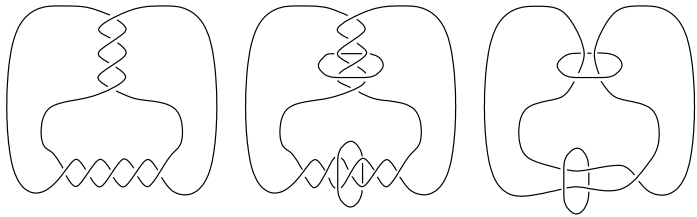}
\end{center}
\caption{Left to right:  The link $K$.  The link $K$ with crossing
circles added.  The link $L$.}
\label{fig:knot-decomp}
\end{figure}

Begin with a prime, twist reduced diagram of the knot $K$.  At each
twist region of the diagram, add a closed curve encircling two
strands, called a \emph{crossing circle}.  We now have a diagram of a
link, which we call $J$.  The link complement $S^3-J$ is homeomorphic
to the complement of the link $L$ obtained from $J$ by removing pairs
of crossings from each twist region.  See Figure
\ref{fig:knot-decomp}.

These links with crossing circles added have been studied by many
people, including Adams \cite{adams:aug}; Lackenby, Agol and Thurston
\cite{lackenby:alt-volume}.  Many of the properties of these links were
addressed in previous papers (see \cite{purcell:volume}, and
\cite{futer-purcell} with Futer).

The links are useful to us in the cone manifold setting because
first, symmetries of the diagram of $L$ allow us to determine the
hyperbolic structure on $S^3-L$ explicitly.  Secondly, $S^3-L$ is
related to the original knot complement via Dehn filling.  In
particular, if $2m_i$ crossings were removed from the $i$-th twist
region to go from the diagram of $J$ to that of $L$, then by
performing a $1/m_i$ Dehn filling on the corresponding $i$-th crossing
circle of $S^3-L$, we re-insert these crossings and put back the
crossing circle, obtaining $S^3-K$ (see e.g.  \cite{rolfsen-book},
Chapter 9).

Thus, starting with a hyperbolic augmented link $L$, we can determine
information including shapes of the cusps.  Next, we can find
conditions which will ensure the Dehn filling above can be obtained by
a hyperbolic Dehn filling via cone deformation.  Hence we will have
the first two pieces needed to bound the normalized height of the cusp
$K$: The initial manifold $M$ with hyperbolic structure will be
$S^3-L$.  It will be connected by cone deformation to $S^3-K$.

\subsubsection{The geometry of $S^3-L$}

We note in \cite{futer-purcell} that if $K$ has at least $2$ twist
regions, then the augmented link complement $S^3-L$ will admit a
complete, finite volume hyperbolic structure.  This can be shown by
work of Adams \cite{adams:aug}, using the geometrization of Haken
manifolds.  It can also be shown more directly using Andreev's
theorem.  The proof using Andreev's theorem gives a packing of circles
related to $S^3-L$ which we will need again in Section
\ref{subsec:2bridge}, so we include it here.

\begin{theorem}
  If $K$ has a prime, twist reduced diagram with at least $2$ twist
  regions, then the augmented link complement $S^3-L$ will admit a
  complete, finite volume hyperbolic structure.
\label{thm:Lhyp}
\end{theorem}

\begin{proof}
  We prove the lemma by finding the hyperbolic structure.  To do so,
  we decompose $S^3-L$ into two ideal polyhedra, as in
  \cite{lackenby:alt-volume} and \cite{futer-purcell}.  We then find a
  hyperbolic structure on the polyhedra which gives a nonsingular
  hyperbolic structure on $S^3-L$ after gluing.

  To find the polyhedra, first remove all remaining single crossings
  at crossing circles.  That is, if the number of crossings in a twist
  region of the diagram of $K$ was odd, then $S^3-L$ will still have a
  single crossing at that twist region.  Remove it.  The new link has
  components lying flat in the projection plane bound together by
  crossing circles.  Call it $L'$.
  
  Slice $S^3-L'$ along the projection plane.  This cuts the manifold
  into two identical pieces, one on either side of the plane.  Each
  2--punctured disk encircled by a crossing circle has been sliced in
  half.  Slice along these remaining halves of disks, opening them up
  into two faces.  This gives an ideal polyhedral decomposition of the
  manifold.  The edges are given by the intersection of the
  2--punctured disks with the projection plane.  On each of the two
  polyhedra, we have one face per planar region of the diagram of
  $L'$, and two triangular faces per 2--punctured disk.  Shade the
  faces arising from 2--punctured disks, giving each polyhedron a
  checkerboard coloring.  The edges meet in 4--valent ideal vertices.
  At each vertex, two white faces are separated by two shaded faces.
  (See \cite{futer-purcell} or \cite{lackenby:alt-volume} for
  pictures.) 
  
  The gluing pattern on the polyhedra is given by following the above
  process in reverse: First, on each polyhedron fold the pairs of
  triangular shaded faces together at their common vertex.  Then glue
  each white face to its corresponding face on the opposite
  polyhedron via the identity map.  This gives back $S^3-L'$.
  
  To obtain $S^3-L$, we use the same two polyhedra, but make the
  following change in gluing.  For crossing circles encircling a
  single crossing of the diagram of $L$, rather than glue shaded
  triangles across their common vertex on a single polyhedron, we
  insert a half--twist by gluing each triangle of one polyhedron to
  the opposite triangle of the other.  See Figure
  \ref{fig:half-twist}.  

\begin{figure}
  \begin{center}
    \includegraphics{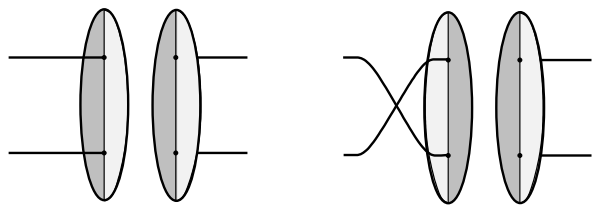}
  \end{center}
  \caption{Left:  Gluing 2--punctured discs with no crossings.  Right:
    Gluing 2--punctured discs with a single crossing.}
  \label{fig:half-twist}
\end{figure}
  
Now our goal is to find a hyperbolic structure on these two polyhedra.
We show the two polyhedra actually have totally geodesic faces in
$\HH^3$, meeting in dihedral angles of $\pi/2$.  Since faces are glued
in pairs, and edges in 4's, the gluing will give a complete hyperbolic
structure on $S^3-L$.
  
Note any totally geodesic plane in $\HH^3$ extends to give a circle on
the boundary $S^2$ at infinity.  Conversely, given a collection of
circles on $S^2$, we can view those circles as boundaries of planes in
$\HH^3$.  We obtain our polyhedra by finding appropriate circle
packings, then cutting away half--spaces bounded by hemispheres in
$\HH^3$.  Our tool is a corollary of Andreev's theorem noted by
Thurston in \cite{thurston}:
  
\begin{theorem}[Andreev]
  Let $\gamma$ be a triangulation of $S^2$ such that each edge has
  distinct ends and no two vertices are joined by more than one edge.
  Then there is a packing of circles in $S^2$ whose nerve is isotopic
  to $\gamma$.  This circle packing is unique up to M{\"o}bius
  transformation.
\end{theorem}

Recall that the \emph{nerve} of a circle packing is the graph obtained
by adding a vertex for each circle, and an edge connecting two
vertices if and only if the corresponding circles are tangent.

We begin by finding a triangulation of $S^2$ associated with one of
the polyhedra described above.  For each white face, take a vertex.
Connect two vertices by an edge if and only if the two corresponding
white faces meet at a vertex of the polyhedron.  This gives a
triangulation $\gamma$ of $S^2$.  Note if we draw vertices of $\gamma$
on top of white faces of the polyhedron, and edges of $\gamma$ through
the vertices of the polyhedron, that each triangle of $\gamma$
circumscribes a shaded face of the polyhedron.

To apply Andreev's theorem, we need to show this triangulation
$\gamma$ satisfies the conditions: each edge has distinct ends, and no
two vertices are joined by more than one edge.  To show these
conditions, we consider the graph dual to $\gamma$, as follows.

Consider the original link diagram.  We view it as a 4--valent graph.
Replace each twist region of the graph by a single edge, as in Figure
\ref{fig:twist-to-edge}.  

\begin{figure}
  \begin{center}
    \includegraphics{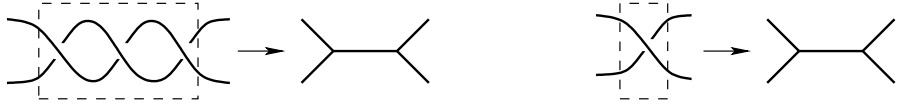}
  \end{center}
  \caption{Create a trivalent graph by replacing each twist region with
  an edge.}
  \label{fig:twist-to-edge}
\end{figure}

Since each crossing is part of some twist region, this gives a
trivalent graph $\Gamma$.
Note the dual graph to $\Gamma$ is isotopic to the triangulation
$\gamma$.
We use properties of the dual $\Gamma$ to show the triangulation
$\gamma$ satisfies the conditions necessary for Andreev's theorem.

\begin{lemma}
  No edge of $\gamma$ has a single vertex at both endpoints.
  \label{lemma:andreev1}
\end{lemma}

\begin{proof}
Consider the graph $\Gamma$.  If its dual $\gamma$ has an edge with one
vertex at both endpoints, then the graph $\Gamma$ must be of the form
in Figure \ref{fig:Andr_double_end} (a).

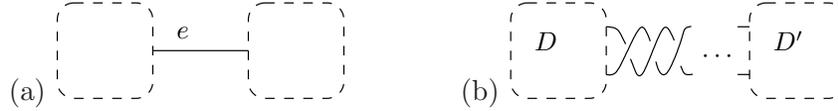
\begin{figure}
  \begin{center}
    (a)
  \input{figures/Andr_double_end.pstex_t}
  \hspace{.5in}
  (b)
  \input{figures/Andr_twist_end.pstex_t}
  \end{center}
  \caption{(a.) Form of $\Gamma$ if $\gamma$ has an edge
    with a single vertex at both endpoints.  (b.)  Corresponding form
    of original link diagram.  }
  \label{fig:Andr_double_end}
\end{figure}

Note then if the edge $e$ was an original edge of the 4--valent link
diagram, then we have a simple closed curve intersecting the link
diagram transversely in a single point.  This is impossible.

If instead the edge $e$ replaced a twist region in the original
4-valent link diagram, then the original diagram had the form of
Figure \ref{fig:Andr_double_end} (b).  However, note one of the
subdiagrams $D$ and $D'$ must contain crossings else the original
diagram has at most one twist region.  Without loss of generality, say
$D$ contains crossings.  Then the dashed simple closed curve enclosing
$D$ intersects the link diagram transversely in two points, but both
its interior and its exterior contain crossings.  This contradicts the
fact that the original diagram was prime.
\end{proof}

\begin{lemma}
  No two vertices of $\gamma$ are joined by more than one edge. 
  \label{lemma:andreev2}
\end{lemma}

\begin{proof}
Suppose two vertices of $\gamma$ are joined by two distinct edges.
Then consider the dual graph $\Gamma$.  $\Gamma$ must be of the form
shown in Figure \ref{fig:Andr_double_edge} (a), where region $B$ and
region $B'$ are the same.

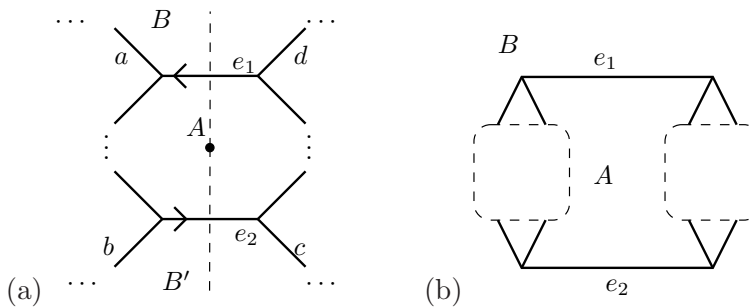
\begin{figure}
  \begin{center}
    (a)
  \input{figures/Andr_double_edge.pstex_t}
  \hspace{.5in}
  (b)
  \input{figures/Andr_double_graph.pstex_t}
 
  \end{center}
  \caption{Form of $\Gamma$ when $\gamma$ has two vertices
  connected by two distinct edges.  In (a), the edges of $\gamma$ are
  the dashed lines, with one of the vertices in $A$ and the other at
  $\infty$.}
  \label{fig:Andr_double_edge}
\end{figure}

Give the edge $e_1$ an orientation, say $e_1$ points in the direction
of the edge $a$.  This gives the region $A$ an orientation, and hence
assigns to the edge $e_2$ an orientation.  Now, trace a path of edges
around the region $B$, listing the edges in the order they occur.  Our
list begins with $e_1$, then $a$, etc.  Note $e_2$ must appear in the
list before the edge $d$, else region $B$ would be cut off from region
$B'$.  Also, edge $b$ must occur before edge $c$, by the orientation
on $e_2$.  Hence the graph $\Gamma$ has the form illustrated in Figure
\ref{fig:Andr_double_edge} (b).

Case 1: Edges $e_1$ and $e_2$ both were edges of the original
4--valent link diagram.  Then there is a simple closed curve
intersecting $e_1$ and $e_2$ transversely in the original link diagram
and intersecting the diagram in no other points.  Since the
original diagram has crossings on either side of the region $A$, this
contradicts the fact that we started with a prime diagram.

Case 2: Edges $e_1$ and $e_2$ both replaced twist regions in the
original link diagram.  Then the original link diagram must have been
of the form in Figure \ref{fig:Andr_double_diag1} (a).  Then since the
original diagram was twist reduced, either the subdiagram $D$ or $D'$
must be a string of bigons.  But then rather than two distinct twist
regions, our diagram must have had only one twist region.  Thus edges
$e_1$ and $e_2$ would not have been distinct.

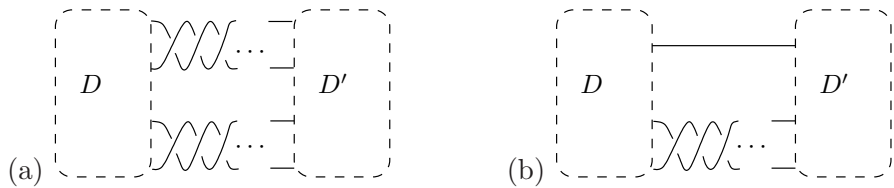
\begin{figure}
  \begin{center}
    (a) \input{figures/Andr_double_diag1.pstex_t}
  \hspace{.5in}
  (b)
  \input{figures/Andr_double_diag2.pstex_t}
  \end{center}
  \caption{Possible forms of the original link diagram if $\gamma$
  contains two vertices connected by two distinct edges.}
  \label{fig:Andr_double_diag1}
\end{figure}

Case 3: One of edge $e_1$ or $e_2$ is from the original link diagram,
the other replaced a twist region.  Then the original link diagram
must have been of the form shown in Figure \ref{fig:Andr_double_diag1}
(b).  Now note that the simple closed curve bounding subdiagram $D$
intersects the link transversely in only three points.  This is
impossible.

Since these are the only three cases, these contradictions prove the
lemma. 
\end{proof}

Now Lemmas \ref{lemma:andreev1} and \ref{lemma:andreev2} are enough to
show that Andreev's theorem applies, and there is a circle packing of
$S^2$ whose nerve is $\gamma$.

Complete the proof of Theorem \ref{thm:Lhyp} by slicing off
half--spaces bounded by geodesic hemispheres in $\mathbb H^3$
corresponding to each circle in the circle packing.  These give the
geodesic white faces of the polyhedron.  The shaded faces are obtained
by slicing off hemispheres in $\mathbb H^3$ corresponding to each
circle of the dual circle packing.  Note white and shaded hemispheres
meet at right angles.  

Thus the polyhedron for the flat augmented link is totally geodesic in
$\mathbb H^3$, and hence the flat augmented link is hyperbolic.
\end{proof}

\subsubsection{Normalized lengths on cusps}  To apply results on cone
deformations, we will need the following information on normalized
lengths of particular curves on the cusps of $S^3-L$.

\begin{prop}
  Let $K$ and $L$ be as above.  In the prime, twist reduced diagram of
  $K$, let $n$ be the number of twist regions.  Let $c_i$ be the
  number of crossings in the $i$-th twist region.  Then on
  cusps of $L$, we have the following normalized lengths.
  \begin{enumerate}
  \item For a cusp corresponding to a crossing circle, let $s_i$ be
    the slope such that {D}ehn filling $S^3-L$ along $s_i$ re-inserts
    the $c_i$ crossings at that twist region.  Then the normalized
    length of $s_i$ is at least $\sqrt{c_i}$.
  \item For the cusp corresponding to the link component in the
    projection plane, the normalized height (i.e. the normalized
    length of the curve orthogonal to a meridian) is at least
    $\sqrt{n}$, and at most $\sqrt{n(n-1)}$.
  \end{enumerate}
\label{prop:proj-height}
\end{prop}

\begin{proof}
  This can be deduced from results of \cite{futer-purcell}.  However,
  given the circle packing associated to $S^3-L$ developed in the
  previous section, we are able to streamline the proof somewhat.
  
  First, recall the polyhedra obtained in the proof of Theorem
  \ref{thm:Lhyp}.  Note each cusp will be tiled by rectangles given by
  the intersection of the cusp with the totally geodesic white and
  shaded faces of the polyhedra.  Two opposite sides of each of these
  rectangles come from the intersection of the cusp with shaded faces
  of the polyhedra, or from the 2--punctured disks in the diagram of
  $L$, and the other two sides come from white faces.  Call the sides
  shaded sides and white sides, respectively.  
  
  For part (1), note any crossing circle intersects a single
  2--punctured disk in a longitude.  Half of the longitude is given by
  the intersection with one polyhedron, the other half by intersection
  with the other.  The crossing circle intersects no other shaded
  faces.  Thus the crossing circle cusp is tiled by exactly two
  identical rectangles.
  
  When the crossing circle encircles no single crossing, or half twist,
  then two white faces intersect the cusp in meridians.  In this case,
  $c_i$ is even, $c_i=2m_i$, and the slope $s_i$ is $1/m_i$.  It is
  given by one step along a white side of a rectangle (a meridian),
  plus $2m_i=c_i$ steps along shaded sides of the rectangle.
  
  When the crossing circle encircles a single crossing, then the two
  rectangles tiling the cusp are sheared.  In this case, a meridian is
  given by one step along a white side, plus (or minus) a step along a
  shaded side.  Now $c_i =2m_i+1$ is odd, and the slope $s_i =1/m_i$
  is given by a step along a white side, plus (or minus) a step along
  a shaded side, plus (or minus) $2m_i$ steps along shaded sides.
  
  In either case, the slope is given by a single step along a white
  side plus (or minus) $c_i$ steps along a shaded side.  Let $w$
  denote the length of a white side, and let $s$ denote the length of
  a shaded side.  Thus the normalized length of $s_i$ is:
  $$\frac{\sqrt{w^2+c_i^2s^2}}{\sqrt{2ws}} =
  \sqrt{\frac{w}{2s}+c_i^2\frac{s}{2w}}.$$
  Let $u=w/(2s)$.  Then for
  constant $c_i$, this quantity is minimized when $u=c_i/2$, and the
  minimum value is $\sqrt{c_i}$.
  
  For part (2), we need to know more about the possible lengths of the
  white sides of rectangles tiling the cusp.  For a given rectangle,
  consider the corresponding vertex in one of the polyhedra.  On the
  sphere $S^2$ at infinity, this vertex will look like a point of
  tangency of two circles in the packing given by Andreev's theorem.
  Since the nerve of this circle packing is a triangulation of $S^2$,
  the point of tangency corresponds to an edge on two distinct
  triangles of this nerve.  Thus there are two additional circles
  tangent to the two given circles, as in Figure \ref{fig:circles}.

\begin{figure}
  \begin{center}
    \includegraphics{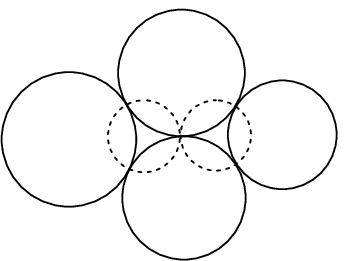} \hspace{.2in}
    \includegraphics{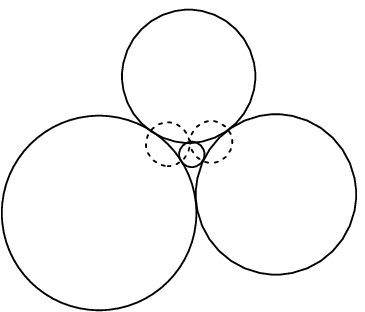}
  \end{center}
  \caption{The form of the packing of circles about any vertex.
    Dashed circles correspond to the shaded faces of the polyhedron at
    this vertex.}
  \label{fig:circles}
\end{figure}

Now, apply a M{\"o}bius transformation taking the vertex, or the point
of tangency of the circle packing, to infinity.  This takes the two
tangent circles to two parallel lines.  It takes the two additional
tangent circles to circles tangent to both the parallel lines, as in
Figure \ref{fig:rec-circles}.  Note this also gives the similarity
structure of the rectangle under consideration.  If we normalize so
that the shaded side (coming from the intersection with a 2--punctured
disk) has length $1$, then the circle lying under the dashed line in
Figure \ref{fig:rec-circles} has diameter $1$.  Since circles in the
circle packing do not overlap, this forces a white side (i.e. without
dashes in the figure) to have length at least $1$.

\begin{figure}
  \begin{center}
    \includegraphics{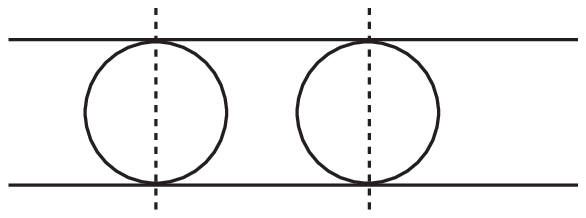}
    \hspace{.2in}
    \includegraphics{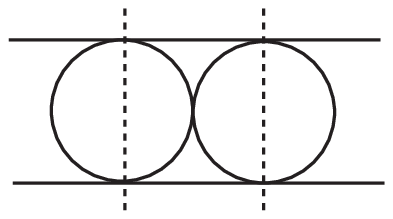}
  \end{center}
  \caption{Result of taking point of tangency to infinity.}
  \label{fig:rec-circles}
\end{figure}

For an upper bound on the length of a white side, note it is longest
when all circles in the circle packing are lined up in a row, tangent
to the two parallel lines.  In this case, when a shaded side has
length $1$, the length of the white side will be the number of circles
between the parallel lines, minus $1$.  To find the number of circles
possible, we use an Euler characteristic argument.  Note there will be
$2n$ shaded faces in one of the polyhedra, coming from opening up $n$
crossing circles.  These each correspond to a single face of the
triangulation of $S^2$ used in Andreev's theorem.  Because it is a
triangulation, the number of edges equals $3/2$ the number of faces,
or $3n$.  Then by an Euler characteristic argument, there are $n+2$
vertices.  Recall each vertex gave a circle in our circle packing.
Since two of these circles become our parallel lines, there are $n$
circles packed between these.  Thus the maximum length of a white side
is $n-1$.

Now, note that the cusp in the projection plane intersects each
2--punctured disk $2n$ times, each time in a meridian.  Half of this
meridian will come from the intersection of one polyhedron (the
``top'' half), the other half from the intersection of the other
polyhedron.  Thus the meridian runs along two shaded sides of
rectangles.  Still normalizing so each shaded side has length $1$,
this implies a meridian has length $2$.  The curve giving the height
$H$ of the cusp runs along $2n$ white sides of the rectangles (one for
each 2--punctured disk of intersection).  The normalized height will
be given by $H/\sqrt{2H} = \sqrt{H/2}.$ The height $H$ will be minimal
when each of the $2n$ white sides it runs along are of minimal length,
or length $1$.  Thus the minimal normalized height is $\sqrt{n}$.

Similarly, the height $H$ will be maximal when each white side is of
maximal length, or length $n-1$.  Thus the maximal normalized height
is $\sqrt{n(n-1)}$.
\end{proof}

\subsubsection{Hyperbolic cone deformation on $S^3-L$}
  \label{subsec:crossings}

If the normalized lengths of the slopes on which {D}ehn filling was
performed were sufficiently large, then Hodgson and Kerckhoff showed
that the {D}ehn filling could be realized by a hyperbolic cone
deformation \cite{hk:univ}.  Proposition \ref{prop:proj-height}
indicates that to attain these minimal normalized slope lengths, we
need only ensure that there are sufficiently many crossings in each
twist region.  In particular, by \cite{hk:univ}, provided $\sqrt{c_j}
\geq 10.6273$, or provided that there are at least $113$ crossings per
twist region of the prime, twist reduced diagram of $K$, then the
{D}ehn filling of $S^3-L$ can be realized by a hyperbolic cone
deformation.

Actually, we will see that the number of crossings in each twist
region needs to be larger than $113$ to apply Theorem
\ref{thm:bound-bdry-time}.  For now, however, note that we may ensure
a hyperbolic cone deformation does exist from $S^3-L$ to $S^3-K$ by
putting conditions on the diagram of $K$.

\subsection{Finding values for tube radius}
\label{subsec:tube}

We need to ensure that the tube radius remains larger than $R_1 =
0.56$ throughout the cone deformation.  We can do so by again applying
work of Hodgson and Kerckhoff \cite{hk:univ}.  They show that in the
neighborhood of the $j$-th component of the singular locus, we can
increase the cone angle from $0$ to $2\pi$, maintaining the tube
radius $R \geq R_1$, provided equation (\ref{eq:Lhat_int}) of the
proof of Lemma \ref{lemma:b_j-bound} holds.

That is, consider a particular component of the singular locus
$\Sigma_j$.  The singular locus consists of crossing circles of the
link $L$.  At the $j$-th crossing circle, the hyperbolic cone
deformation performs {D}ehn filling along a slope of the form $1/n_j$.
Let $L_0$ denote the normalized length of this slope in $S^3-L$.  We
can then increase the cone angle from $0$ to $2\pi$, maintaining the
tube radius bound $R \geq R_1$, provided the normalized length $L_0^2$
is larger than the increasing function $I(R_1)$ defined in equation
(\ref{eq:Lhat_int}).

So to ensure $R \geq R_1 \geq 0.56$, we compute that $L_0^2$ must be
at least $113.044$.  By Proposition \ref{prop:proj-height}, we know
$L_0^2$ is at least $c_i$, where $c_i$ is the number of crossings in
the $i$-th twist region.  Thus by ensuring there are at least $114$
crossings in the $i$-th twist region, we will guarantee that our tube
radius is at least $0.56$.  For general $R_1$, we have the following
result.

\begin{lemma}
  Let $K$, $L$ be as above, such that we have a cone deformation from
  $S^3-L$ to $S^3-K$.  Fix $R_1 \geq 0.531$.
  Then $R \geq R_1$ throughout the deformation provided the number of
  crossings in each twist region is at least $I(R_1)$, where $I(R_1)$
  is defined in equation (\ref{eq:Lhat_int}).
\label{lemma:crossings-vs-R}
\end{lemma}

\subsection{Proof of Theorem \ref{thm:main-knot-cusp}}
\label{subsec:numeric}

\hspace{.1in}

\medskip

\noindent\textbf{Theorem \ref{thm:main-knot-cusp}} \emph{
  Let $K$ be a knot in $S^3$ which admits a prime, twist reduced
  diagram with a total of $n\geq 2$ twist regions, with each twist
  region containing at least $c\geq 116$ crossings.  In a hyperbolic
  structure on $S^3-K$, take a horoball neighborhood $U$ about $K$.
  Normalize so that the meridian on $\partial U$ has length $1$.  Then
  the height $H$ of the cusp of $S^3-K$ satisfies:
  $$H\geq n\,(1-f(c))^2.$$
  Here $f(c)$ is a positive function of $c$ which
  approaches $0$ as $c$ increases to infinity.
Additionally, $H\leq n(\sqrt{n-1} + f(c))^2$.
}

\begin{proof}
  
  If $K$ admits a prime, twist reduced diagram with at least $2$ twist
  regions, then we may form the augmented link $L$ from the diagram of
  $K$.  We have seen in \S\ref{subsec:crossings} that, provided there
  are at least $113$ crossings in each twist region of the diagram of
  $K$, then $S^3-L$ and $S^3-K$ are connected by a hyperbolic cone
  deformation $X_t$, with $X_0=S^3-L$, $X_\tau=S^3-K$.  The singular
  locus will consist of the crossing circle components of $L$.  There
  will be $n$ of these.
  
  By results of Section \ref{subsec:tube}, if there are at least $114$
  crossings per twist region then the tube radius of the deformation
  is at least $R_1=0.56$.  Then by Theorem \ref{thm:main-norm-height}, 
  $$ -(2\pi)^2\frac{\sqrt{n\,C(R_1)}}{(1-e^{-2R_1})\sqrt{2}}
  \leq h_\tau-h_0
  \leq
  (2\pi)^2\frac{\sqrt{n\,C(R_1)}}{(1-e^{-2R_1})\sqrt{2}},$$
  where $h_t$ denotes the normalized height of the cusp at time $t$.
  
And so:
\begin{equation}
 h_0 - (2\pi)^2 \frac{\sqrt{n\,C(R_1)}}{(1-e^{-2R_1})\sqrt{2}} \leq
h_\tau \leq
h_0 + (2\pi)^2 \frac{\sqrt{n\,C(R_1)}}{(1-e^{-2R_1})\sqrt{2}}.
\label{eqn:2-sides}
\end{equation}

We may now apply part (2) of Proposition \ref{prop:proj-height}.  This
implies that the initial normalized height $h_0$, i.e. the normalized
height on $S^3-L$, is between $\sqrt{n}$ and $\sqrt{n(n-1)}$.

Plug this into the left hand side of the above equation and we find:
\begin{equation}
h_\tau \geq \sqrt{n}\left(1-(2\pi)^2
  \frac{\sqrt{C(R_1)}}{(1-e^{-2R_1})\sqrt{2}}\right).
\label{eqn:h}
\end{equation}

For the right hand side, we find:
\begin{equation}
  h_\tau \leq \sqrt{n}\left(\sqrt{n-1} + (2\pi)^2
  \frac{\sqrt{C(R_1)}}{(1-e^{-2R_1})\sqrt{2}}\right).
\label{eqn:hright}
\end{equation}

The term $C(R_1)$ can be computed explicitly using the formulas in the
proof of 
Lemma \ref{lemma:b_j-bound}.  Recall $C(R_1)$ is strictly decreasing,
approaching $0$ as $R_1$ approaches infinity.  Thus the function
$$\bar{f}(R_1)= (2\pi)^2\frac{\sqrt{C(R_1)}}{(1-e^{-2R_1})\sqrt{2}},$$
which appears in both (\ref{eqn:h}) and (\ref{eqn:hright}), is
decreasing with $R_1$, approaching $0$.  The graph of $\bar{f}$, as
well as that of $I(R)$, is shown in Figure \ref{fig:F-I-graphs}.

\begin{figure}
  \begin{center}
    \framebox{
      \begin{tabular}{c}
    \includegraphics[height=1.6in]{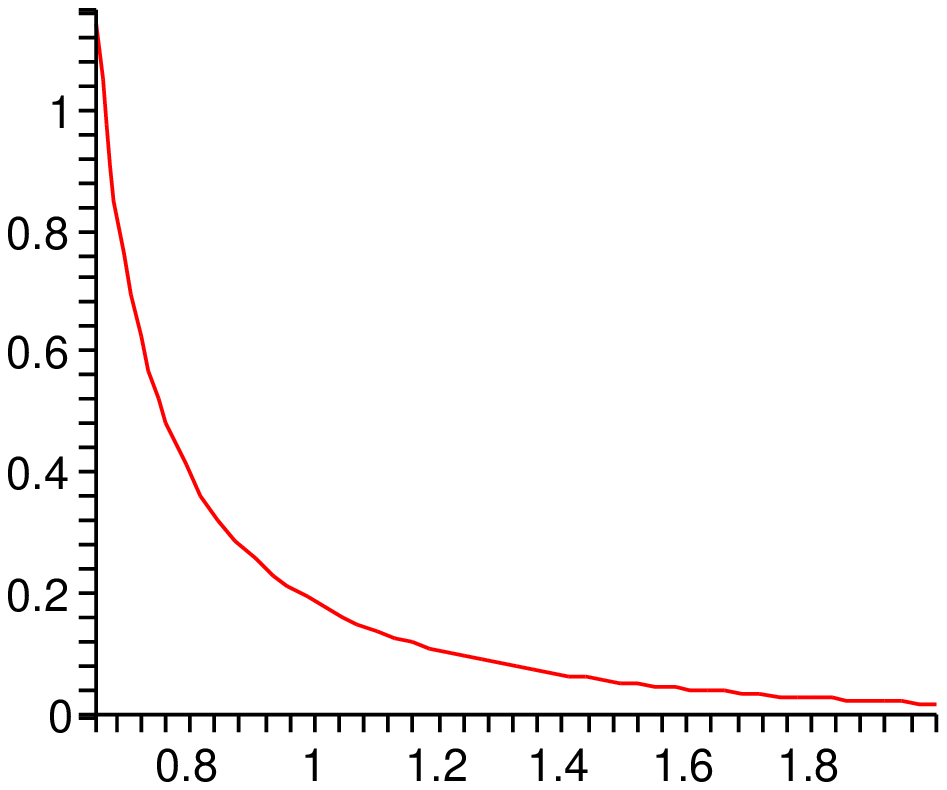}\\
    $y=\bar{f}(R)$
    \end{tabular}}
\hspace{.1in}
\framebox{
  \begin{tabular}{c}
    \includegraphics[height=1.6in]{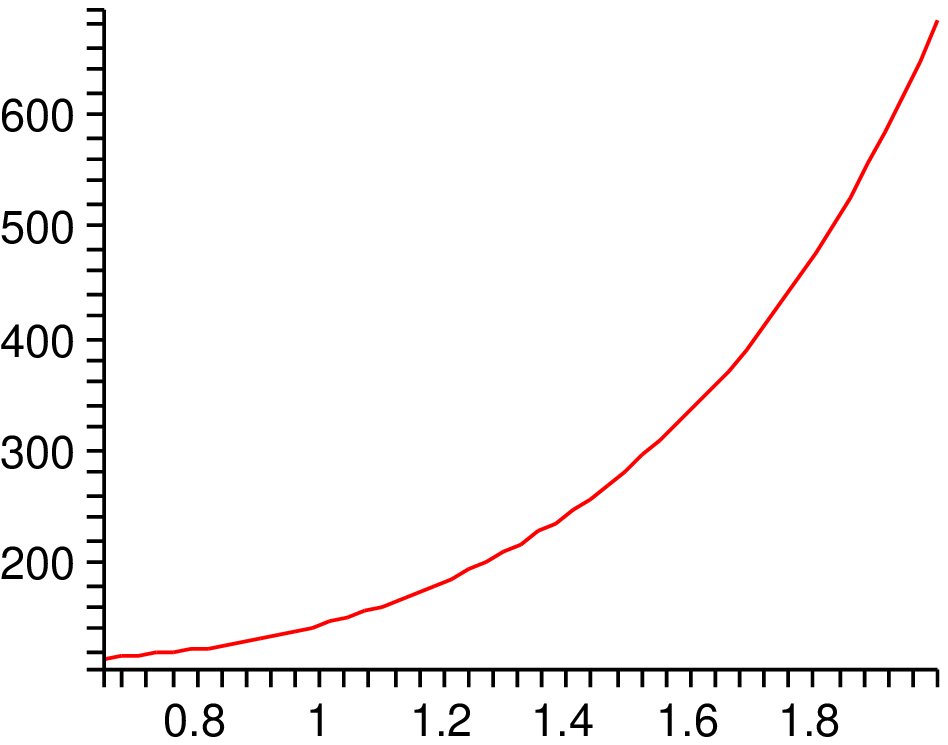}\\
    $y=I(R)$
  \end{tabular}}
  \end{center}
  \caption{Graphs of $\bar{f}(R)$ and $I(R)$.  Here $\bar{f}(R)$ gives the
    bounds in Theorem \ref{thm:main-knot-cusp} for given tube radius
    $R$.  $I(R)$ gives the number of crossings required to guarantee
    tube radius $R$.}
  \label{fig:F-I-graphs}
\end{figure}

When $R_1=0.56$, $1-\bar{f}(R_1)$ is negative.  Since $h_\tau$
is known to be positive, this doesn't give a very good estimate in
(\ref{eqn:h}).  To get a better estimate, we increase $R_1$.  The
value $1-\bar{f}(R_1)$ becomes positive when $R_1$ is about $0.6624$.
However, as $R_1$ increases, so does $I(R_1)$, so by Lemma
\ref{lemma:crossings-vs-R}, we need to increase the number of
crossings in each twist region to guarantee that the tube radius
remains larger than $R_1$ throughout the deformation.  In particular,
$c$ must be at least $116$ to ensure $R_1 \geq 0.6624$.

In general, let $f(c) = \bar{f}(I^{-1}(c))$.  Since $I(R)$ is
strictly increasing for $R\geq 0.56$, this is well defined.  Then for
$c\geq 116$, $1-f(c)$ is positive.

To finish the proof, note that we can scale so that the meridian of
the cusp of $S^3-K$ has length $1$.  Then the actual height $H$ of the
cusp of $S^3-K$ will be at least $(h_\tau)^2$.
\end{proof}

Notice in Figure \ref{fig:F-I-graphs} that $\bar{f}$ is close to $0$
for values of $R$ near $2$.  However for $R$ this size, $I(R)$, or the
number of crossings required per twist region, is nearly $700$.  We
find a nice compromise by choosing $R=1.0$.  Then $I(R)$ is not much
larger than $116$, and yet $1-\bar{f}$ is already over $0.8$.  We use
this fact in the following corollary, which gives more explicit
numbers to the results of Theorem \ref{thm:main-knot-cusp}.

\begin{corollary}
  Let $K$ be as in the previous theorem, with at least $c\geq 145$
  crossings in each twist region.  Then when the meridian of the cusp
  of $S^3-K$ is normalized to have length $1$, the height of the cusp
  satisfies: $$H \geq n\,(0.81544)^2.$$
\label{cor:0.8}
\end{corollary}

\begin{proof}
  Letting $c\geq 145$ in the above theorem ensures that $R> R_1 =
  1.0$.  Then $1-f(c) = 1-\bar{f}(1.0) \geq 0.81544$.
\end{proof}

\subsection{Example: 2--bridge knots}
\label{subsec:2bridge}

The requirement of Theorem \ref{thm:main-knot-cusp} that there be at
least $116$ crossings per twist region was based on the worst case
analysis.  This came from Proposition \ref{prop:proj-height}, in which
we showed that in the worst possible case, the normalized length of a
slope was at least the square root of the number of crossings.  When
we are dealing with particular classes of knots, we can often make
significantly better estimates, and reduce the number of crossings
required.  In this section, we present a particular example of such a
class of knots:  2--bridge knots.

\begin{figure}
  \begin{center}
  \input{figures/2bridge.pstex_t}
  \end{center}
  \caption{Augmented link of a 2--bridge knot with an even number
    of twist regions, and an even number of crossings per region.}
  \label{fig:2bridge}
\end{figure}
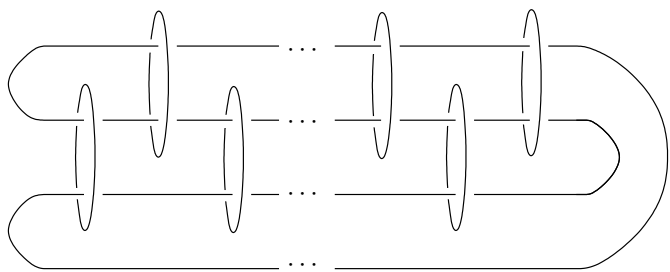

We will determine the circle packing obtained from the proof of
Theorem \ref{thm:Lhyp} for a 2--bridge knot.  A 2--bridge knot with an
even number of twist regions will have a corresponding augmented link
$L$ of the form of Figure \ref{fig:2bridge}, except possibly with
single crossings added at some crossing circles.  The picture for an
odd number of twist regions is similar.

\begin{figure}
  \begin{center}
  \input{figures/2bridge_tri.pstex_t}
  \end{center}
  \caption{Trivalent graph associated with a 2--bridge knot.}
  \label{fig:2bridge_tri}
\end{figure}
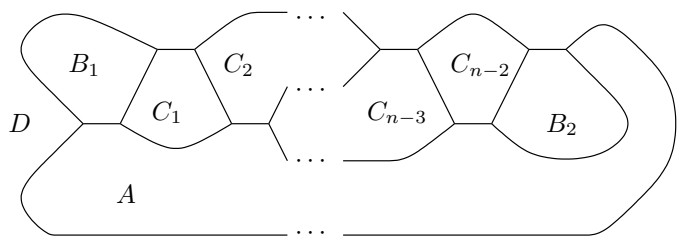

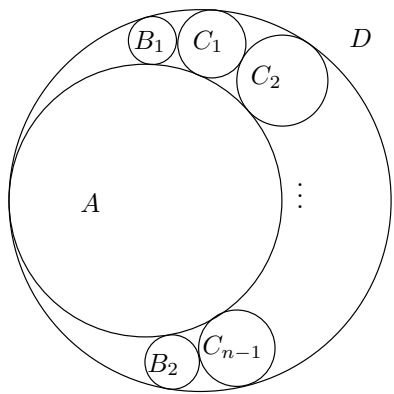
\begin{figure}
  \begin{center}
  \input{figures/Fig_2bridge_circ.pstex_t}
  \end{center}
  \caption{Circle packing associated with a 2--bridge knot.}
  \label{fig:2bridge_circ}
\end{figure}

In any case, the link will correspond to a trivalent graph $\Gamma$ of
the form shown in Figure \ref{fig:2bridge_tri}.  Note each region of
the graph complement gives a circle in the circle packing, with
tangencies across edges. Thus each circle will be tangent to the
circles corresponding to regions $A$ and $D$.  The circle packing
associated with this graph is of the form in Figure
\ref{fig:2bridge_circ}.

Recall that the shape of cusps can be determined by analyzing the
rectangles we get from taking points of tangency of the circles in
Figure \ref{fig:2bridge_circ} to infinity.  In particular, we are
interested in the cusp shapes which correspond to crossing circles in
the link $L$.  The rectangles tiling these cusps come from the
tangencies of circle $A$ with certain of the $B_i$ and $C_i$, and from
tangencies of $D$ with certain of the $B_i$ and $C_i$.  Thus we find
the general shape of rectangles of this form.

When we take a point of tangency between $A$ and $B_i$ to infinity,
$i=1$ or $2$, the two circles tangent to both $A$ and $B_i$ are $D$
and a $C_j$.  These are tangent to each other.  Thus our rectangle in
this case is a square.  Similarly for points of tangency between $D$
and $B_i$.

A crossing circle coming from one of these points of tangency will
therefore be tiled by two squares, with a longitude running across two
shaded sides.  Then the slope $s_i$ along which Dehn filling is
performed will have normalized length $\sqrt{1+c_i^2}/\sqrt{2}$.  Note
this is larger than $I(0.6624)$, implying that the results of Theorem
\ref{thm:main-knot-cusp} apply, provided $c_i \geq 16$.  It is larger
than $I(1.0)$, implying the results of Corollary \ref{cor:0.8} hold,
when $c_i \geq 17$.  Thus the number of crossings in our theorem has
gone from $116$ to $16$, and $145$ to $17$.

When we take a point of tangency between $A$ and $C_i$, or between $D$
and $C_i$ to infinity, we get a $1$ by $2$ rectangle.  This is because
the circles $A$ and $C_i$ are tangent to two distinct circles
$C_{i-1}$ and $C_{i+1}$, or possibly a $B_i$, which are each tangent
to $D$.  Since $D$ is also tangent to $A$ and $C_i$, in the picture at
infinity we have two parallel lines, $A$ and $C_i$, with three full
sized circles tangent to each other and also to the lines $A$ and
$C_i$.

In this case, the crossing circle is tiled by two $1$ by $2$
rectangles, stacked to form a square.  The slope $s_i$
will have normalized length $\sqrt{4+c_i^2}/\sqrt{4}$.  It will be
larger than $I(0.6624)$ provided $c_i \geq 22$, and larger
than $I(1.0)$ provided $c_i \geq 24$.

We can also make improvements on the initial normalized height,
i.e. the value $h(X)$ of Theorem \ref{thm:main-norm-height}.  To
compute $h(X)$, we determine the shape of the rectangles tiling the
cusp corresponding to the knot strand.  These include all rectangles
obtained by taking a point of tangency of the circle packing of Figure
\ref{fig:2bridge_circ} to infinity, except those tangencies
corresponding to crossing circle cusps.

Each circle tangency of the form $C_{2k}$ tangent to $D$ will give a
rectangle on the knot strand cusp.  Similarly, those tangencies of
circles $C_{2k+1}$ tangent to $A$ will give rectangles on the knot
strand cusp.  As we saw above, these rectangles will each have
dimensions $1$ by $2$.  The side of length $2$ will contribute to the
height.  There will be $(n-2)$ of these rectangles, since each $C_j$
contributes one.  Thus these rectangles contribute $2n-4$ to the
height.

Each circle tangency $C_i$ tangent to $C_{i+1}$ also gives a rectangle
on the knot strand cusp, as do the tangencies $B_1$ with $C_1$ and
$C_{n-2}$ with $B_2$.  Each of these rectangles is a square,
contributing $1$ to the height, and there are $(n-1)$ such rectangles.

The tangencies of $B_1$ and $D$ and of $B_2$ and $A$ (or $D$ if $n$ is
odd) also contribute rectangles of dimension $1$ by $1$ to the knot
strand.  Finally, the tangency of $A$ to $D$ gives a rectangle with
dimensions $1$ by $(n-1)$, contributing $(n-1)$ to the height.  Thus
the total height of the cusp is
$$(2n-4) + (n-1) + 2 + (n-1) = 4n -4.$$

(When $n=2$ this argument needs to be modified slightly:  there are no
$C_i$'s, but $B_1$ and $B_2$ are tangent.  Tangencies of $A$ and
$B_2$, $D$ and $B_1$, and $A$ and $D$ also contribute to the height.
These all give $1$ by $1$ rectangles, hence the height is $4 = 4n-4$
in this case as well.)

A meridian runs along the widths of two rectangles in the cusp
tiling.  Hence it has length $2$.  Thus the normalized height $h(X)$
is given by:
$$h(X) = \frac{4n-4}{\sqrt{(4n-4)(2)}} = \sqrt{2(n-1)}.$$

This information will show:

\begin{prop}
  Let $K$ be a 2--bridge knot in $S^3$ which admits a prime, twist
  reduced diagram with $n\geq 2$ twist regions and at least $c\geq 24$
  crossings in each twist region.  In a hyperbolic structure on
  $S^3-K$, normalized such that the meridian has length $1$, the
  height of the cusp of $S^3-K$ satisfies:
  $$\left(\sqrt{2(n-1)} - \sqrt{n}(0.18456)\right)^2 \leq H \leq
	\left(\sqrt{2(n-1)} + \sqrt{n}(0.18456)\right)^2.$$
\label{prop:2bridge}
\end{prop}

\begin{proof}
Let $L$ be the augmented link corresponding to $K$.  Provided we have
at least $n\geq 2$ twist regions, Theorem \ref{thm:Lhyp} implies
$S^3-L$ is hyperbolic.  We have seen in the above discussion that,
provided we have at least $c\geq 24$ crossings in each twist region,
the slope on each crossing circle of $L$ on which we perform Dehn
filling has normalized length at least $I(1.0)$.  This implies that a
cone deformation exists from $S^3-L$ to $S^3-K$ for which the tube
radius is at least $R_1= 1.0$ throughout.  Then Theorem
\ref{thm:main-norm-height} implies the normalized height of the cusp
of $S^3-K$ satisfies:
$$-\frac{\sqrt{n}(2\pi)^2\sqrt{C(1.0)}}{(1-e^{-2.0})} \leq
 h_{\tau} - \sqrt{2(n-1)} \leq
 \frac{\sqrt{n}(2\pi)^2\sqrt{C(1.0)}}{(1-e^{-2.0})}.$$
Using the formulas for $C(1.0)$ from the proof of Lemma
 \ref{lemma:b_j-bound}, we find
$$\sqrt{2(n-1)} - \sqrt{n}(0.18456) \leq h_{\tau} \leq
 \sqrt{2(n-1)} + \sqrt{n}(0.18456).$$
We obtain the final result by noting if we scale such that the
meridian has length $1$, then the height is at least $(h_\tau)^2$.
\end{proof}

%% file: figures/Andr_double_end.pstex_t
\begin{picture}(0,0)%
\includegraphics{figures/Andr_double_end.pstex}%
\end{picture}%
\setlength{\unitlength}{3947sp}%
\begingroup\makeatletter\ifx\SetFigFont\undefined%
\gdef\SetFigFont#1#2#3#4#5{%
  \reset@font\fontsize{#1}{#2pt}%
  \fontfamily{#3}\fontseries{#4}\fontshape{#5}%
  \selectfont}%
\fi\endgroup%
\begin{picture}(1824,624)(589,-373)
\put(1351, 14){\makebox(0,0)[lb]{\smash{{\SetFigFont{10}{12.0}{\familydefault}{\mddefault}{\updefault}{\color[rgb]{0,0,0}$e$}%
}}}}
\end{picture}%

%% file: figures/Andr_twist_end.pstex_t
\begin{picture}(0,0)%
\includegraphics{figures/Andr_twist_end.pstex}%
\end{picture}%
\setlength{\unitlength}{3947sp}%
\begingroup\makeatletter\ifx\SetFigFont\undefined%
\gdef\SetFigFont#1#2#3#4#5{%
  \reset@font\fontsize{#1}{#2pt}%
  \fontfamily{#3}\fontseries{#4}\fontshape{#5}%
  \selectfont}%
\fi\endgroup%
\begin{picture}(2124,624)(589,-373)
\put(751,-61){\makebox(0,0)[lb]{\smash{{\SetFigFont{10}{12.0}{\familydefault}{\mddefault}{\updefault}{\color[rgb]{0,0,0}$D$}%
}}}}
\put(2251,-61){\makebox(0,0)[lb]{\smash{{\SetFigFont{10}{12.0}{\familydefault}{\mddefault}{\updefault}{\color[rgb]{0,0,0}$D'$}%
}}}}
\put(1801,-136){\makebox(0,0)[lb]{\smash{{\SetFigFont{10}{12.0}{\familydefault}{\mddefault}{\updefault}{\color[rgb]{0,0,0}$\cdots$}%
}}}}
\end{picture}%

%% file: figures/Andr_double_edge.pstex_t
\begin{picture}(0,0)%
\includegraphics{figures/Andr_double_edge.pstex}%
\end{picture}%
\setlength{\unitlength}{3947sp}%
\begingroup\makeatletter\ifx\SetFigFont\undefined%
\gdef\SetFigFont#1#2#3#4#5{%
  \reset@font\fontsize{#1}{#2pt}%
  \fontfamily{#3}\fontseries{#4}\fontshape{#5}%
  \selectfont}%
\fi\endgroup%
\begin{picture}(1612,1828)(211,-1466)
\put(601, 14){\makebox(0,0)[lb]{\smash{{\SetFigFont{10}{12.0}{\familydefault}{\mddefault}{\updefault}{\color[rgb]{0,0,0}$a$}%
}}}}
\put(826,239){\makebox(0,0)[lb]{\smash{{\SetFigFont{10}{12.0}{\familydefault}{\mddefault}{\updefault}{\color[rgb]{0,0,0}$B$}%
}}}}
\put(226,239){\makebox(0,0)[lb]{\smash{{\SetFigFont{10}{12.0}{\familydefault}{\mddefault}{\updefault}{\color[rgb]{0,0,0}$\cdots$}%
}}}}
\put(1726, 14){\makebox(0,0)[lb]{\smash{{\SetFigFont{10}{12.0}{\familydefault}{\mddefault}{\updefault}{\color[rgb]{0,0,0}$d$}%
}}}}
\put(1801,239){\makebox(0,0)[lb]{\smash{{\SetFigFont{10}{12.0}{\familydefault}{\mddefault}{\updefault}{\color[rgb]{0,0,0}$\cdots$}%
}}}}
\put(1801,-586){\makebox(0,0)[lb]{\smash{{\SetFigFont{10}{12.0}{\familydefault}{\mddefault}{\updefault}{\color[rgb]{0,0,0}$\vdots$}%
}}}}
\put(1051,-436){\makebox(0,0)[lb]{\smash{{\SetFigFont{10}{12.0}{\familydefault}{\mddefault}{\updefault}{\color[rgb]{0,0,0}$A$}%
}}}}
\put(526,-586){\makebox(0,0)[lb]{\smash{{\SetFigFont{10}{12.0}{\familydefault}{\mddefault}{\updefault}{\color[rgb]{0,0,0}$\vdots$}%
}}}}
\put(526,-1186){\makebox(0,0)[lb]{\smash{{\SetFigFont{10}{12.0}{\familydefault}{\mddefault}{\updefault}{\color[rgb]{0,0,0}$b$}%
}}}}
\put(301,-1411){\makebox(0,0)[lb]{\smash{{\SetFigFont{10}{12.0}{\familydefault}{\mddefault}{\updefault}{\color[rgb]{0,0,0}$\cdots$}%
}}}}
\put(901,-1411){\makebox(0,0)[lb]{\smash{{\SetFigFont{10}{12.0}{\familydefault}{\mddefault}{\updefault}{\color[rgb]{0,0,0}$B'$}%
}}}}
\put(1726,-1186){\makebox(0,0)[lb]{\smash{{\SetFigFont{10}{12.0}{\familydefault}{\mddefault}{\updefault}{\color[rgb]{0,0,0}$c$}%
}}}}
\put(1801,-1411){\makebox(0,0)[lb]{\smash{{\SetFigFont{10}{12.0}{\familydefault}{\mddefault}{\updefault}{\color[rgb]{0,0,0}$\cdots$}%
}}}}
\put(1341,-18){\makebox(0,0)[lb]{\smash{{\SetFigFont{10}{12.0}{\familydefault}{\mddefault}{\updefault}{\color[rgb]{0,0,0}$e_1$}%
}}}}
\put(1354,-1098){\makebox(0,0)[lb]{\smash{{\SetFigFont{10}{12.0}{\familydefault}{\mddefault}{\updefault}{\color[rgb]{0,0,0}$e_2$}%
}}}}
\end{picture}%

%% file: figures/Andr_double_graph.pstex_t
\begin{picture}(0,0)%
\includegraphics{figures/Andr_double_graph.pstex}%
\end{picture}%
\setlength{\unitlength}{3947sp}%
\begingroup\makeatletter\ifx\SetFigFont\undefined%
\gdef\SetFigFont#1#2#3#4#5{%
  \reset@font\fontsize{#1}{#2pt}%
  \fontfamily{#3}\fontseries{#4}\fontshape{#5}%
  \selectfont}%
\fi\endgroup%
\begin{picture}(1824,1657)(1189,-1460)
\put(1351, 89){\makebox(0,0)[lb]{\smash{\SetFigFont{10}{12.0}{\familydefault}{\mddefault}{\updefault}{\color[rgb]{0,0,0}$B$}%
}}}
\put(1951, 14){\makebox(0,0)[lb]{\smash{\SetFigFont{10}{12.0}{\familydefault}{\mddefault}{\updefault}{\color[rgb]{0,0,0}$e_1$}%
}}}
\put(1951,-736){\makebox(0,0)[lb]{\smash{\SetFigFont{10}{12.0}{\familydefault}{\mddefault}{\updefault}{\color[rgb]{0,0,0}$A$}%
}}}
\put(2026,-1411){\makebox(0,0)[lb]{\smash{\SetFigFont{10}{12.0}{\familydefault}{\mddefault}{\updefault}{\color[rgb]{0,0,0}$e_2$}%
}}}
\end{picture}

%% file: figures/Andr_double_diag1.pstex_t
\begin{picture}(0,0)%
\includegraphics{figures/Andr_double_diag1.pstex}%
\end{picture}%
\setlength{\unitlength}{3947sp}%
\begingroup\makeatletter\ifx\SetFigFont\undefined%
\gdef\SetFigFont#1#2#3#4#5{%
  \reset@font\fontsize{#1}{#2pt}%
  \fontfamily{#3}\fontseries{#4}\fontshape{#5}%
  \selectfont}%
\fi\endgroup%
\begin{picture}(2124,1074)(1189,-1573)
\put(1351,-1036){\makebox(0,0)[lb]{\smash{{\SetFigFont{10}{12.0}{\familydefault}{\mddefault}{\updefault}{\color[rgb]{0,0,0}$D$}%
}}}}
\put(2851,-1036){\makebox(0,0)[lb]{\smash{{\SetFigFont{10}{12.0}{\familydefault}{\mddefault}{\updefault}{\color[rgb]{0,0,0}$D'$}%
}}}}
\put(2326,-811){\makebox(0,0)[lb]{\smash{{\SetFigFont{10}{12.0}{\familydefault}{\mddefault}{\updefault}{\color[rgb]{0,0,0}$\cdots$}%
}}}}
\put(2326,-1411){\makebox(0,0)[lb]{\smash{{\SetFigFont{10}{12.0}{\familydefault}{\mddefault}{\updefault}{\color[rgb]{0,0,0}$\cdots$}%
}}}}
\end{picture}%

%% file: figures/Andr_double_diag2.pstex_t
\begin{picture}(0,0)%
\includegraphics{figures/Andr_double_diag2.pstex}%
\end{picture}%
\setlength{\unitlength}{3947sp}%
\begingroup\makeatletter\ifx\SetFigFont\undefined%
\gdef\SetFigFont#1#2#3#4#5{%
  \reset@font\fontsize{#1}{#2pt}%
  \fontfamily{#3}\fontseries{#4}\fontshape{#5}%
  \selectfont}%
\fi\endgroup%
\begin{picture}(2124,1074)(1189,-1573)
\put(1351,-1036){\makebox(0,0)[lb]{\smash{{\SetFigFont{10}{12.0}{\familydefault}{\mddefault}{\updefault}{\color[rgb]{0,0,0}$D$}%
}}}}
\put(2851,-1036){\makebox(0,0)[lb]{\smash{{\SetFigFont{10}{12.0}{\familydefault}{\mddefault}{\updefault}{\color[rgb]{0,0,0}$D'$}%
}}}}
\put(2326,-1411){\makebox(0,0)[lb]{\smash{{\SetFigFont{10}{12.0}{\familydefault}{\mddefault}{\updefault}{\color[rgb]{0,0,0}$\cdots$}%
}}}}
\end{picture}%

%% file: figures/2bridge.pstex_t
\begin{picture}(0,0)%
\includegraphics{figures/2bridge.pstex}%
\end{picture}%
\setlength{\unitlength}{3947sp}%
\begingroup\makeatletter\ifx\SetFigFont\undefined%
\gdef\SetFigFont#1#2#3#4#5{%
  \reset@font\fontsize{#1}{#2pt}%
  \fontfamily{#3}\fontseries{#4}\fontshape{#5}%
  \selectfont}%
\fi\endgroup%
\begin{picture}(4164,1712)(1239,-1766)
\put(3001,-1711){\makebox(0,0)[lb]{\smash{{\SetFigFont{10}{12.0}{\familydefault}{\mddefault}{\updefault}{\color[rgb]{0,0,0}$\cdots$}%
}}}}
\put(3001,-1261){\makebox(0,0)[lb]{\smash{{\SetFigFont{10}{12.0}{\familydefault}{\mddefault}{\updefault}{\color[rgb]{0,0,0}$\cdots$}%
}}}}
\put(3001,-811){\makebox(0,0)[lb]{\smash{{\SetFigFont{10}{12.0}{\familydefault}{\mddefault}{\updefault}{\color[rgb]{0,0,0}$\cdots$}%
}}}}
\put(3001,-361){\makebox(0,0)[lb]{\smash{{\SetFigFont{10}{12.0}{\familydefault}{\mddefault}{\updefault}{\color[rgb]{0,0,0}$\cdots$}%
}}}}
\end{picture}%

%% file: figures/2bridge_tri.pstex_t
\begin{picture}(0,0)%
\includegraphics{figures/2bridge_tri.pstex}%
\end{picture}%
\setlength{\unitlength}{3947sp}%
\begingroup\makeatletter\ifx\SetFigFont\undefined%
\gdef\SetFigFont#1#2#3#4#5{%
  \reset@font\fontsize{#1}{#2pt}%
  \fontfamily{#3}\fontseries{#4}\fontshape{#5}%
  \selectfont}%
\fi\endgroup%
\begin{picture}(4216,1528)(1186,-1841)
\put(1201,-1111){\makebox(0,0)[lb]{\smash{{\SetFigFont{10}{12.0}{\familydefault}{\mddefault}{\updefault}{\color[rgb]{0,0,0}$D$}%
}}}}
\put(1576,-736){\makebox(0,0)[lb]{\smash{{\SetFigFont{10}{12.0}{\familydefault}{\mddefault}{\updefault}{\color[rgb]{0,0,0}$B_1$}%
}}}}
\put(2101,-1036){\makebox(0,0)[lb]{\smash{{\SetFigFont{10}{12.0}{\familydefault}{\mddefault}{\updefault}{\color[rgb]{0,0,0}$C_1$}%
}}}}
\put(2551,-736){\makebox(0,0)[lb]{\smash{{\SetFigFont{10}{12.0}{\familydefault}{\mddefault}{\updefault}{\color[rgb]{0,0,0}$C_2$}%
}}}}
\put(3001,-436){\makebox(0,0)[lb]{\smash{{\SetFigFont{10}{12.0}{\familydefault}{\mddefault}{\updefault}{\color[rgb]{0,0,0}$\cdots$}%
}}}}
\put(3001,-886){\makebox(0,0)[lb]{\smash{{\SetFigFont{10}{12.0}{\familydefault}{\mddefault}{\updefault}{\color[rgb]{0,0,0}$\cdots$}%
}}}}
\put(3001,-1336){\makebox(0,0)[lb]{\smash{{\SetFigFont{10}{12.0}{\familydefault}{\mddefault}{\updefault}{\color[rgb]{0,0,0}$\cdots$}%
}}}}
\put(3001,-1786){\makebox(0,0)[lb]{\smash{{\SetFigFont{10}{12.0}{\familydefault}{\mddefault}{\updefault}{\color[rgb]{0,0,0}$\cdots$}%
}}}}
\put(1876,-1561){\makebox(0,0)[lb]{\smash{{\SetFigFont{10}{12.0}{\familydefault}{\mddefault}{\updefault}{\color[rgb]{0,0,0}$A$}%
}}}}
\put(3976,-736){\makebox(0,0)[lb]{\smash{{\SetFigFont{10}{12.0}{\familydefault}{\mddefault}{\updefault}{\color[rgb]{0,0,0}$C_{n-2}$}%
}}}}
\put(4576,-1111){\makebox(0,0)[lb]{\smash{{\SetFigFont{10}{12.0}{\familydefault}{\mddefault}{\updefault}{\color[rgb]{0,0,0}$B_2$}%
}}}}
\put(3451,-1036){\makebox(0,0)[lb]{\smash{{\SetFigFont{10}{12.0}{\familydefault}{\mddefault}{\updefault}{\color[rgb]{0,0,0}$C_{n-3}$}%
}}}}
\end{picture}%

%% file: figures/Fig_2bridge_circ.pstex_t
\begin{picture}(0,0)%
\includegraphics{figures/Fig_2bridge_circ.pstex}%
\end{picture}%
\setlength{\unitlength}{3947sp}%
\begingroup\makeatletter\ifx\SetFigFont\undefined%
\gdef\SetFigFont#1#2#3#4#5{%
  \reset@font\fontsize{#1}{#2pt}%
  \fontfamily{#3}\fontseries{#4}\fontshape{#5}%
  \selectfont}%
\fi\endgroup%
\begin{picture}(2416,2414)(1193,-1868)
\put(3341,304){\makebox(0,0)[lb]{\smash{{\SetFigFont{10}{12.0}{\familydefault}{\mddefault}{\updefault}{\color[rgb]{0,0,0}$D$}%
}}}}
\put(1651,-736){\makebox(0,0)[lb]{\smash{{\SetFigFont{10}{12.0}{\familydefault}{\mddefault}{\updefault}{\color[rgb]{0,0,0}$A$}%
}}}}
\put(1986,294){\makebox(0,0)[lb]{\smash{{\SetFigFont{10}{12.0}{\familydefault}{\mddefault}{\updefault}{\color[rgb]{0,0,0}$B_1$}%
}}}}
\put(2356,294){\makebox(0,0)[lb]{\smash{{\SetFigFont{10}{12.0}{\familydefault}{\mddefault}{\updefault}{\color[rgb]{0,0,0}$C_1$}%
}}}}
\put(2721, 69){\makebox(0,0)[lb]{\smash{{\SetFigFont{10}{12.0}{\familydefault}{\mddefault}{\updefault}{\color[rgb]{0,0,0}$C_2$}%
}}}}
\put(3006,-706){\makebox(0,0)[lb]{\smash{{\SetFigFont{10}{12.0}{\familydefault}{\mddefault}{\updefault}{\color[rgb]{0,0,0}$\vdots$}%
}}}}
\put(2416,-1626){\makebox(0,0)[lb]{\smash{{\SetFigFont{10}{12.0}{\familydefault}{\mddefault}{\updefault}{\color[rgb]{0,0,0}$C_{n-1}$}%
}}}}
\put(2071,-1736){\makebox(0,0)[lb]{\smash{{\SetFigFont{10}{12.0}{\familydefault}{\mddefault}{\updefault}{\color[rgb]{0,0,0}$B_2$}%
}}}}
\end{picture}%

%% file: jpurcell-cusps.bbl
\begin{thebibliography}{10}

\bibitem{adams:aug}
Colin~C. Adams, \emph{Augmented alternating link complements are hyperbolic},
  Low-dimensional topology and Kleinian groups (Coventry/Durham, 1984), London
  Math. Soc. Lecture Note Ser., vol. 112, Cambridge Univ. Press, Cambridge,
  1986, pp.~115--130, MR0903861, Zbl 0632.57008.

\bibitem{adams:waist}
\bysame, \emph{Waist size for cusps in hyperbolic 3-manifolds}, Topology
  \textbf{41} (2002), no.~2, 257--270, MR1876890, Zbl 0985.57012.

\bibitem{agol:bounds}
Ian Agol, \emph{Bounds on exceptional {D}ehn filling}, Geom. Topol. \textbf{4}
  (2000), 431--449, MR1799796, Zbl 0959.57009.

\bibitem{ast:volumes}
Ian Agol, Peter~A. Storm, and William~P. Thurston, \emph{{Lower bounds on
  volumes of hyperbolic Haken 3-manifolds}}, J. Amer. Math. Soc. \textbf{20}
  (2007), no.~4, 1053--1077, with an appendix by Nathan
Dunfield, MR2328715, Zbl pre05177853.

\bibitem{boroczky}
K.~B{\"o}r{\"o}czky, \emph{Packing of spheres in spaces of constant curvature},
  Acta Math. Acad. Sci. Hungar. \textbf{32} (1978), no.~3-4, 243--261,
	MR0512399, Zbl 0422.52011.

\bibitem{brock-bromberg-evans-souto:tame}
Jeffrey Brock, Kenneth Bromberg, Richard Evans, and Juan Souto, \emph{Tameness
  on the boundary and {A}hlfors' measure conjecture}, Publ. Math. Inst. Hautes
  \'Etudes Sci. (2003), no.~98, 145--166, MR2031201, Zbl 1060.30054.

\bibitem{brock-bromberg:density}
Jeffrey~F. Brock and Kenneth~W. Bromberg, \emph{On the density of geometrically
  finite {K}leinian groups}, Acta Math. \textbf{192} (2004), no.~1,
33--93, MR2079598, Zbl 1055.57020.

\bibitem{bromberg:cone}
K.~Bromberg, \emph{Hyperbolic cone-manifolds, short geodesics, and {S}chwarzian
  derivatives}, J. Amer. Math. Soc. \textbf{17} (2004), no.~4, 783--826
  (electronic), MR2083468, Zbl 1061.30037.

\bibitem{bromberg:cone-rigid}
\bysame, \emph{Rigidity of geometrically finite hyperbolic cone-manifolds},
  Geom. Dedicata \textbf{105} (2004), 143--170, MR2057249, Zbl 1057.53029.

\bibitem{eudave-munoz-luecke}
Mario Eudave-Mu{\~n}oz and J.~Luecke, \emph{Knots with bounded cusp volume yet
  large tunnel number}, J. Knot Theory Ramifications \textbf{8} (1999), no.~4,
  437--446, MR1697382, Zbl 0942.57005.

\bibitem{fkp:volumes}
David Futer, Efstratia Kalfagianni, and Jessica~S. Purcell,
\emph{{Dehn filling, volume, and the Jones polynomial}},
J. Differential Geom. \textbf{78} (2008), no.~3, 429--464, MR2396249,
Zbl pre05261987. 

\bibitem{futer-purcell}
David Futer and Jessica~S. Purcell, \emph{{Links with no exceptional
  surgeries}}, Comment. Math. Helv. \textbf{82} (2007), no.~3,
629--628, MR2314056, Zbl 1134.57003.

\bibitem{hk:cone-rigid}
Craig~D. Hodgson and Steven~P. Kerckhoff, \emph{Rigidity of hyperbolic
  cone-manifolds and hyperbolic {D}ehn surgery}, J. Differential Geom.
  \textbf{48} (1998), no.~1, 1--59, MR1622600, Zbl 0919.57009.

\bibitem{hk:univ}
\bysame, \emph{Universal bounds for hyperbolic {D}ehn surgery}, Ann. of Math.
  (2) \textbf{162} (2005), no.~1, 367--421, MR2178964, Zbl 1087.57011.

\bibitem{lackenby:surg}
Marc Lackenby, \emph{Word hyperbolic {D}ehn surgery}, Invent. Math.
  \textbf{140} (2000), no.~2, 243--282, MR1756996, Zbl 0947.57016.

\bibitem{lackenby:alt-volume}
\bysame, \emph{The volume of hyperbolic alternating link complements}, Proc.
  London Math. Soc. (3) \textbf{88} (2004), no.~1, 204--224, With an appendix
  by Ian Agol and Dylan Thurston, MR2018964, Zbl 1041.57002.

\bibitem{lehto}
Olli Lehto, \emph{Univalent functions and {T}eichm\"uller spaces}, Graduate
  Texts in Mathematics, vol. 109, Springer-Verlag, New York, 1987,
	MR0867407, Zbl 0606.30001.

\bibitem{lickorish}
W.~B.~R. Lickorish, \emph{A representation of orientable combinatorial
  {$3$}-manifolds}, Ann. of Math. (2) \textbf{76} (1962), 531--540,
MR0151948, Zbl 0106.37102.

\bibitem{mostow}
G.~D. Mostow, \emph{Strong rigidity of locally symmetric spaces}, Princeton
  University Press, Princeton, N.J., 1973, Annals of Mathematics Studies, No.
  78, MR0385004, Zbl 0265.53039.

\bibitem{prasad}
Gopal Prasad, \emph{Strong rigidity of {${\bf Q}$}-rank {$1$} lattices},
  Invent. Math. \textbf{21} (1973), 255--286, MR0385005, Zbl 0264.22009.

\bibitem{purcell:volume}
Jessica~S. Purcell, \emph{{Volumes of highly twisted knots and links}}, Algebr.
  Geom. Topol. \textbf{7} (2007), 93--108, MR2289805, Zbl 1135.57005.

\bibitem{rolfsen-book}
Dale Rolfsen, \emph{Knots and links}, Publish or Perish Inc., Berkeley, Calif.,
  1976, Mathematics Lecture Series, No. 7, MR0515288, Zbl 0339.55004.

\bibitem{thurston}
William~P. Thurston, \emph{The geometry and topology of three-manifolds},
  Princeton Univ. Math. Dept. Notes, 1979.

\bibitem{thurston:bulletin}
\bysame, \emph{Three-dimensional manifolds, {K}leinian groups and hyperbolic
  geometry}, Bull. Amer. Math. Soc. (N.S.) \textbf{6} (1982), no.~3,
357--381, MR648524, Zbl 0528.57009.

\bibitem{wallace}
Andrew~H. Wallace, \emph{Modifications and cobounding manifolds}, Canad. J.
  Math. \textbf{12} (1960), 503--528, MR0125588, Zbl 0108.36101.

\end{thebibliography}
